\documentclass[twocolumn]{autart}    % Enable this line and disable the 
\usepackage{graphicx}          % Include this line if your 
\usepackage{float}             % document contains figures,
 \usepackage[caption=false]{subfig}

\usepackage{amsmath,amsfonts,epsfig,array}
\usepackage{cite}
\usepackage[outline]{contour}
\usepackage{color}
\usepackage[normalem]{ulem}
\newcommand{\chgk}[1]{ {\color{red} #1}}

\usepackage{cite}
\usepackage{comment}
\begin{document}

\begin{frontmatter}

\title{ Boundary Control and Observer Design Via Backstepping for a Coupled Parabolic-Elliptic System} % Title, preferably not more 

\author[Waterloo]{Ala' Alalabi}\ead{aalalabi@uwaterloo.ca},    % Add the 
\author[Waterloo]{Kirsten Morris}\ead{kmorris@uwaterloo.ca}              % e-mail address 
%\author[Baiae]{Publius Maro Vergilius}\ead{vergilius@culture.ir}  % (ead) as shown

\address[Waterloo]{Department of Applied Mathematics,\\
        University of Waterloo, 200 University Avenue West, \\
        Waterloo, ON, Canada  N2L 3G1}  % Please supply                                              
\thanks{ *This work was supported by NSERC and the Faculty of Mathematics, University of Waterloo}% <-this % stops a space

                             % keyword list or with the 
                                          % help of the Automatica 
                                          % keyword wizard

\begin{abstract}                          % Abstract of not more than 200 words.

Stabilization  of a coupled system consisting of a parabolic partial differential equation and an elliptic partial differential equation is considered. Even in the situation when the parabolic equation is exponentially stable on its own,  the coupling between the two equations can cause instability in the overall system.  A backstepping approach is used to derive a boundary control input that  stabilizes the coupled system. The result is an explicit expression for the stabilizing control law. The second part  of the paper involves  the design of  exponentially convergent observers to estimate the state of the coupled system,  given some partial boundary measurements. The observation error system is shown to be exponentially stable, again  by employing a backstepping method. This leads to  the design of observer gains in closed-form. Finally, we address the output-feedback problem by combining the  observers  with the state feedback boundary control. The theoretical results are demonstrated with numerical simulations.
\end{abstract}
\end{frontmatter}

%%%%%%%%%%%%%%%%%%%%%%%%%%%%%%%%%%%%%%%%%%%%%%%%%%%%%%%%%%%%%%%%%%%%%%%%%%%%%%%%
\section{Introduction}

%%%% Paragraph including real-life appliactions of parabolic-elliptic systems + well-posedness of such systems
The coupling of parabolic partial differential equations with elliptic partial differential equations appears in  a number of physical applications, including  electro-chemical models of lithium-ion  cells \cite{consiglieri2020weak,wu2006well}, biological transport networks \cite{li2020global},  chemotaxis phenomena \cite{ahn2019global,negreanu2019parabolic} and the thermistor \cite{meinlschmidt2017optimal}. Parabolic-elliptic systems  are thus  an important class of  partial differential-algebraic equations (PDAEs).  Well-posedness of linear PDAEs has been addressed in the literature \cite{trostorff2020semigroups,thaller1996factorization,sviridyuk2003linear,jacob2022solvability}. 
There have also been research efforts dedicated to investigate the well-posdeness of particular linear parabolic-elliptic systems \cite{ahn2019global,benilan1996mild}. In addition,  the well-posedness of quasilinear and nonlinear parabolic-elliptic systems have received a recent attention; see  \cite{wang2016quasilinear,negreanu2019parabolic,coclite2005wellposedness,malysheva2020global,li2020global}. 

%%%%%%% (1) state the overall problem
Parabolic-elliptic systems can exhibit instability, leading to several complications  in real-life applications. One example can be found in the context of chemotaxis phenomena. In this case the parabolic equation  describes the diffusion of cells while the elliptic equation models the  concentration of  a chemical attractant. In \cite{tao2015boundedness} authors showed that the aforementioned model can manifest  unstable  dynamics, leading to spatially complex patterns in cell density as well as  the concentration of the chemical; see \cite{polezhaev2003spatial}. These patterns will in turn impact many biological processes including  tissue development, tumor growth and wound healing. A more recent illustration for the instability in coupled parabolic-elliptic systems can be found in \cite{ParadaCerpaMorris}. The findings in this paper demonstrated that even when the parabolic equation is exponentially stable, coupling with the elliptic equation  can lead  to an unstable system. 

%%%%%%% (2) review the literature 
Stabilization through boundary control for  coupled linear parabolic partial differential equations has been studied in the literature. In \cite{baccoli2015boundary}, the backstepping approach was used to stabilize the dynamics of a linear coupled reaction-diffusion systems with constant coefficients. An extension of this work to systems with variable coefficients was presented in \cite{Vzquez2016BoundaryCO}. Koga et al. \cite{koga2016backstepping} described boundary control of the one phase Stefan problem,  modeled  by a diffusion equation coupled with an ordinary differential equation. Feedback stabilization of a PDE-ODE system was also studied in \cite{hasan2016boundary}.  

  \begin{sloppypar}
There is less research  tackling  the stabilization   of coupled parabolic-elliptic systems. \cite{krstic2008boundary,krstic2006backstepping,ParadaCerpaMorris}. In \cite[chap. 10]{krstic2008boundary}, stabilization of  parabolic-elliptic systems arose in the context  of stabilizing  boundary control of linearized Kuramoto-Sivashinsky and  Korteweg de Vries equations. The controller  required the presence of two  Dirichlet control inputs. A similar approach was also followed in \cite{krstic2006backstepping} where a hyperbolic-elliptic system arose in the control  of a  Timoshenko beam.  More recently, Parada et al. \cite{ParadaCerpaMorris} considered the boundary control of  unstable parabolic-elliptic system through Dirichlet control with input delay.  
%%%%% (3) outline the contributions of the paper briefly but clearly. 
  \end{sloppypar}
%%%% First contribution 
The first main contribution of this paper is  design of a single feedback Neumann control input that  exponentially stabilizes the dynamics of the two coupled equations.  %Alternatively, it  will enhance the decay rate of the solutions when the system is already stable.
The control input is directly designed on the system of partial differential equations without  approximation by finite-dimensional systems. This will be done by  a backstepping approach \cite{krstic2008boundary}. Backstepping is one of the few methods that yields an explicit control law for PDEs without first approximating the PDE.  These transformations are generally formulated as a Volterra operator, which guarantees under weak conditions the invertibility of the transformation.  When using backstepping,  it is typical to determine the destabilizing terms in the system, find a suitable exponentially stable target system where the destabilizing terms are eliminated   by state transformation and feedback control.  Then, look for an invertible state transformation of the original  system into the exponentially stable target system. This requires finding a kernel of the Volterra operator and also showing that the kernel is well-defined as the solution of an auxilary PDE. One possible approach for stabilization of a parabolic-elliptic system is  to convert  the coupled system into one equation in terms of the parabolic state. However, this will result in the presence of a Fredholm operator  that makes it difficult to establish a suitable kernel for the backstepping transformation.  Another approach would be a vector-valued transformation for both of the parabolic and the elliptic states but this is quite complex. A different approach is taken here.  We  use  a single  transformation previously used for a parabolic equation \cite{krstic2008boundary}. Properties of the kernel of this transformation have already been established.  This  leads to an  unusual target system  in a parabolic-elliptic form.  As in other backstepping designs,   an explicit expression for the controller is obtained as a byproduct of the transformation. Then the  problem is  to  establish the stability of the  target system obtained  from the transformation, which will imply the stability of the original coupled system via the invertible transformation.
%The problem is then to either establishing the stability of the  target system obtained  from the transformation, which will imply the stability of the original coupled system; or to show that the controlled original system is exponentially stable. 
Explicit calculation of the eigenfunctions is not required.

In many dynamical systems, the full state is not  available. This issue  motivated the study of constructing an estimate of the state by designing an observer. The literature on  observer  design for coupled systems appears to focus on the observer
synthesis of systems governed by coupled hyperbolic PDEs or coupled parabolic PDEs; see for instance \cite{aamo2012disturbance,moura2013observer,vazquez2011backstepping,baccoli2015anticollocated}.  
There have been  few  papers addressing  observer design problem for partial differential equations coupled with an elliptic equation \cite{jiang2019two,krstic2006backstepping}. In \cite{jiang2019two} authors designed a
state observer for a coupled  parabolic-elliptic system by requiring a two-sided boundary input for the observer. Also, in the same work mentioned earlier for  Timoshenko beam  \cite{krstic2006backstepping}, authors studied its observer design  by using a hyperbolic-elliptic system and requiring two control inputs. In \cite{smyshlyaev2005backstepping} an observer is designed for a parabolic equation that includes a Volterra term. 
We design  an observer using two measurements, and also several observers that require only one measurement.
 The exponential stability of the  observation error dynamics is achieved by means of designing suitable output injections, also known as filters or observer gains.  In parallel to our approach for the boundary controller design, we derive observer gains 
by using  backstepping transformations that are well-established in the literature. The transformations result in target error systems whose stability is studied. %or that the original observation error dynamics with the obtained  output injections is stable.  

We finally combine the state feedback and observer designs to obtain an  output feedback controller for the  coupled parabolic-elliptic system. Numerical simulations are presented that  illustrate the theoretical findings for each of the objectives.
To our knowledge, this paper is the first to  designing a state estimator for a coupled parabolic-elliptic system with partial boundary measurements,  and only one control input. 

This paper is structured as follows. Section \ref{section1} presents the well-posedness of the parabolic-elliptic systems under consideration. Stability analysis for the uncontrolled system is also described.  Section \ref{section2} includes  the first main result which is the use of a backstepping method to design a boundary controller for the coupled
system. A preliminary version of section \ref{section2}, on stabilization, is included in the proceedings of the 2023 Conference on Decision and Control \cite{cdc}. The design of a state observer for the coupled system is given in Section \ref{section3}. Different designs are proposed based on the available measurements. The output feedback problem is described in Section \ref{section4}. Conclusions and some discussion of possible extensions are given in Section \ref{section5}.

%%%%%%%%%%%%%%%%%%%%%%%%%%%%%%%%%%%%%%%%%%%%%%%%%%%%%%%%%%%%%%%%%%%%%%%%%%%%%%%%%%%%
\section{Well-posedness and stability of system} \label{section1}

We study parabolic-elliptic systems of the form 
\begin{align}
    w_t(x,t) =&  w_{xx}(x,t) - \rho w(x,t) + \alpha v (x,t), \label{1}\\
   0 =&  v_{xx}(x,t) - \gamma v(x,t) + \beta w(x,t), \label{2}\\
   w_x(0,t)=& 0, \quad w _x(1,t) =u(t), \label{3}\\
    v_x (0,t)=&0, \quad v_x (1,t) =0,  \label{4}
\end{align}
where $x \in [0,1]$ and $t \geq 0$. The parameters $\rho, \; \alpha, \; \beta, \; \gamma $ are all real, with $\alpha$, $\beta$ both nonzero. 
%The given restriction on $\gamma$ ensures the well-posedness of system  \eqref{1}-\eqref{4}; see \cite{jacob2022solvability}. 

 With the notation 
$\Delta v(x)=\frac{d^2 v}{dx} ,$ define   the operator $A^{\gamma} : D( A^{\gamma} ) \to L^2(0,1)$
\begin{align*}
   & A^{\gamma} v(x)= (\gamma I -\Delta) v(x) = w(x), \\
    &D( A^{\gamma} ) = \{ v \in H^2 (0,1) , \; v^\prime (0) = v^\prime (1) = 0 \}.
\end{align*}
For values of $\gamma$ that are not eigenvalues of  $A^\gamma$, that is   $\gamma \neq -n \pi^2 $ with $n=0, \dots$, the inverse operator of $A^{\gamma}$   $(A^{\gamma})^{-1}: L^2(0,1) \to D( A^{\gamma} )  $ exists. 
In this situation, the uncontrolled system \eqref{1}-\eqref{4} is well-posed \cite{jacob2022solvability}. Alternatively, define
\begin{align}
   & A= \Delta-\rho I + \alpha \beta (\gamma I - \Delta)^{-1} , \label{10} \\ 
    & D(A) = \{ w \in H^2 (0,1) , \; w^\prime (0) = w^\prime (1) = 0  \}.\nonumber
\end{align}

\begin{thm} \label{C0 generator}
If $\gamma \neq -(n \pi)^2$ the  operator $A$ generates a $C_0$- semigroup and  the control system \eqref{1}-\eqref{4} with observation $w(0,t)$ is well-posed on the state-space $L^2 (0,1).$ It is similarly well-posed with control instead at $x=0$, and/or observation at $x=1 .$
\end{thm}
\noindent \textbf{Proof.} With $\alpha=0$ the control system is the  heat equation with Neumann boundary control. This control system is well-known to be well-posed on $L^2(0,1)$ . Since the operator $A$ is a bounded perturbation of $\Delta$, then the conclusion of the theorem  follows. $\hfill \square $

It will henceforth be assumed that  $\gamma \neq - ( n \pi)^2.$

\begin{thm} \label{th1}
Let $u(t) \equiv 0$. The eigenvalues of system  \eqref{1}-\eqref{4} operator $A$  are
\begin{align}
    \lambda_n = -\rho + \frac{\alpha \beta}{\gamma + (n\pi)^2} - (n \pi)^2, \; \;  n=0, \; 1, \; \dots \label{eigenvalues}
\end{align}

\end{thm}
\noindent \textbf{Proof.} The analysis  is standard but given for completeness. Let $\{\phi_j\}_{j \geq 0} \subset \mathcal{C}^4(0,1)$ be the eigenfunctions of the operator $A$  corresponding to the eigenvalues $\lambda_j$, then setting  $\beta (\gamma I - \partial_{xx})^{-1} \phi_j = e_j$,
\begin{align}
     & \lambda_j \phi_j(x) =  \phi_j^{''}(x)  -\rho  \phi_j(x) + \alpha e_j(x)  \label{11} \\
     & 0 = e_j^{''}(x)  -\gamma  e_j(x) + \beta \phi_j(x) \label{12} \\
     & \phi_j^{'} (0)=\phi_j^{'} (1)=0\\
     & e_j^{'} (0)=e_j^{'} (1)=0.
\end{align}
Solving \eqref{11} for $e_j(x)$
\begin{align}
     e_j (x) = \frac{\rho + \lambda_j}{\alpha \beta } \phi_j (x) - \frac{1}{\alpha \beta} \phi_j^{''}(x) . \label{e_j} 
\end{align}
Substituting for $e_j(x)$ in \eqref{12}, we obtain the  fourth-order differential equation
\begin{align}
     \phi_j^{''''}(x)  & -  (\lambda_j + \rho +\gamma) \phi_j^{''}(x) \nonumber \\
     & + (\gamma (\lambda_j + \rho)- \alpha \beta)\phi_j (x) = 0, \label{13}
\end{align}
with the  boundary conditions
\begin{align}
     \phi_j^{'} (0) =  \phi_j^{'} (1) =\phi_j^{'''} (0) =  \phi_j^{'''} (1)=0. \label{14} 
\end{align}
Solving system \eqref{13}-\eqref{14} for $\phi_j$ yields that $\phi_j= \cos( j \pi x)$ for $j=0,\;1,\dots $. Subbing $\phi_j$ in \eqref{13} and solving for $\lambda_j$ leads to 
\eqref{eigenvalues}. \hfill   $\square$

\begin{cor} \label{Th2}
System (\ref{1})-~(\ref{4}) is exponentially stable if and only if 
\begin{align}
   \rho > \frac{\alpha \beta}{\gamma}, \label{16}
\end{align}
and the decay rate in that case is bounded by the maximum eigenvalue 
$ \rho - \frac{\alpha \beta}{\gamma}$.
\end{cor}
\noindent \textbf{Proof.} 
Since $\Delta $ with domain $D(A)$ is a a Riesz-spectral operator, then  since $A$  is a bounded perturbation, it is also a spectral operator. Alternatively, we note that $A$ is a self-adjoint operator with a compact inverse and hence it is Riesz-spectral \cite[section 3]{curtain2020introduction}. Thus,  $A$ generates a $C_0$-semigroup with growth  bound  determined by the eigenvalues.
$\hfill   \square$

Thus, even in the case when the parabolic equation is exponentially stable, coupling with the elliptic system  can cause the uncontrolled system  to be unstable.

\section{Stabilization } \label{section2}

To design a stabilizing control input, a backstepping approach will be used. Unlike the work done in \cite[Chap.10]{krstic2008boundary} where two control inputs are used, we stabilize the dynamics of the coupled parabolic-elliptic equations \eqref{1}-\eqref{4} using  a single control signal.  The following lemma will be used.
\begin{lem} \cite[chap. 4]{krstic2008boundary} \label{Wellposedness}
    For any $c_2 >0$  the hyperbolic partial differential equation
    \begin{subequations} \label{26}
\begin{align}
    & k^a_{yy}(x,y) - k^a_{xx}(x,y) + c_2 k^a(x,y) = 0,  \quad 0<y<x<1 \label{26a}\\
    & k^a_y(x,0)=0, \quad  k^a(x,x)= -\frac{1}{2} c_2 x, \label{26b}
\end{align}
\end{subequations}
has a continuous unique solution   
\begin{align}
    k^a(x,y)= -c_2 x \frac{I_1 \left( \sqrt{c_2 (x^2-y^2)} \right)}{\sqrt{c_2 (x^2-y^2)}}, \label{k-kernel}
\end{align}
 where $I_1(\cdot)$ is the modified Bessel function of first order defined as \begin{align*}
    I_1(x)= \sum_{m=0}^{\infty}  \frac{(x/2)^{2m+1}}{m! (m+1)!}.
\end{align*}
\end{lem}
We  apply the invertible state transformation
 \begin{align}
        \tilde{w}(x,t)=&w(x,t)-\int_0^x k^a(x,y) w(y,t)dy, \label{-18}
\end{align}
on the parabolic state $w(x,t)$, while the elliptic  state $v(x,t)$ is unchanged.  Here the kernel of the transformation $k^a(x,y)$ is given by \eqref{k-kernel}. The inverse transformation of \eqref{-18} was given in \cite{krstic2008boundary}.
 \begin{lem} \cite[chap. 4]{krstic2008boundary}
     The inverse transformation of (\ref{-18}) is 
     \begin{align}
        w(x,t)=&\tilde{w}(x,t)+\int_0^x \ell^a (x,y) \tilde{w}(y,t)dy,  \label{inverse of transformation 1}
     \end{align}
where $l(x,y)$ is the solution of the  system
\begin{subequations} \label{inverse systemab}
 \begin{align}
     &\ell^a_{xx}(x,y)-\ell^a_{yy}(x,y) +c_2 \ell^a(x,y) =0, \label{inverse system-a} \\
     & \ell^a_y(x,0)=0, \quad \ell^a(x,x)=-\frac{1}{2} c_2 x, \label{inverse system-b}
 \end{align}
 \end{subequations}
 that is
 \begin{align}
    \ell^a(x,y)= -c_2 x \frac{J_1 \left( \sqrt{c_2 (x^2-y^2)} \right)}{\sqrt{c_2 (x^2-y^2)}}, \label{28}
\end{align}
where  $J_1(\cdot)$ is the Bessel function of first order defined as
  \begin{align*}
    J_1(x)= \sum_{m=0}^{\infty}  (-1)^m\frac{(x/2)^{2m+1}}{m! (m+1)!}.
\end{align*}
 \end{lem}   

In what follows we set $c_2=c_1-\rho$ with $c_2>0$.
\begin{thm} \label{thm1}
If the control  signal $u(t)$ is given by 
\begin{align}
        u(t)=&\int_0^1 k^a_x(1,y) w(y,t)dy+k^a(1,1) w(1,t),   \label{27}
\end{align}
    then transformation \eqref{-18}, with $k^a(x,y)$ given by system \eqref{26},  converts the parabolic-elliptic system  \eqref{1}-\eqref{4} into the  target system
\begin{align}
     \tilde{w}_t(x,t) =&  \tilde{w}_{xx}(x,t) - (c_2+\rho) \tilde{w}(x,t)  + \alpha v(x,t)  \nonumber \\ 
     & - \alpha \int_0^x k^a(x,y) v(y,t) dy, \label{19}\\
    0 =&  v_{xx}(x,t)  - \gamma v(x,t) + \beta \tilde{w}(x,t)   \nonumber \\ 
     & + \beta  \int_0^x \ell^a (x,y) \tilde{w}(y,t)  dy, \label{20}\\
   \tilde{w}_x(0,t) =&0, \quad \tilde{w}_x(1,t) =0, \label{21}\\
   v_x (0,t) =&0, \quad v_x (1,t) =0.  \label{22}
\end{align}
\end{thm}
\noindent \textbf{Proof.}
It will prove useful to rewrite \eqref{-18} as
 \begin{align}
      w(x,t)  =&\tilde{w}(x,t)+\int_0^x k^a(x,y) w(y,t)dy. \label{18}
\end{align}
We differentiate \eqref{18} with respect to $x$ twice
%\begin{align}
       %& w_{x}(x,t) =   \tilde{w}_{x}(x,t) + \int_0^x k^a_{x}(x,y)w(y,t)dy  +  k^a(x,x) \nonumber \\
        %& \times w(x,t). \hfill \label{section2-first derivative}
%\end{align}
\begin{align}
       & w_{xx}(x,t) =  \tilde{w}_{xx}(x,t) + \int_0^x k^a_{xx}(x,y)w(y,t)dy \nonumber \\ 
     &+k^a_x(x,x)  w(x,t) + \frac{d}{dx} k^a(x,x) w(x,t) \nonumber \\
     & + k^a(x,x) w_x(x,t), \hfill \label{23}
\end{align}
and with respect to $t$
\begin{align}
         &w_{t}(x,t)=\tilde{w}_{t}(x,t)+\int_0^x k^a(x,y) w_t(y,t)dy \nonumber \\
    % &= \tilde{w}_{t}(x,t) + \int_0^x k^a(x,y) [w_{yy}(y,t) - \rho w(y,t)  + \alpha v(y,t)]dy \nonumber \\
      &=   \tilde{w}_{t}(x,t) + k^a(x,x) w_x(x,t)  - \int_0^x k^a_y(x,y) w_{y}(y,t) dy  \nonumber \\
      &  - \rho \int_0^x k^a(x,y) w(y,t) dy + \alpha \int_0^x k^a(x,y) v(y,t)dy \nonumber \\
         &=    \tilde{w}_{t}(x,t) + k^a(x,x) w_x(x,t) - k^a_y(x,x) w(x,t)   \nonumber \\
         & +k^a_y(x,0)w(0,t) +\int_0^x k^a_{yy}(x,y) w(y,t) dy \nonumber \\
         & -\rho \int_0^x k^a(x,y) w(y,t)dy  + \alpha \int_0^x k^a(x,y) v(y,t)dy. 
         \label{24}
\end{align}

Here 
\begin{align*}
& k^a_x(x,x) = \frac{\partial}{\partial x} k^a(x,y) |_{x=y}, \;  k^a_y(x,x) = \frac{\partial}{\partial y} k^a(x,y) |_{x=y},\\
&\frac{d}{dx}k^a(x,x)= k^a_x(x,x) +k^a_y(x,x). 
\end{align*}
Substituting ~(\ref{23}) and ~(\ref{24}) in ~(\ref{1}), 
\begin{align}
      & \tilde{w}_{t}(x,t) + k^a(x,x) w_x(x,t) - k^a_y(x,x) w(x,t)   \nonumber \\
         & +k^a_y(x,0)w(0,t) +\int_0^x k^a_{yy}(x,y) w(y,t) dy  \nonumber \\ 
     & -\rho \int_0^x k^a(x,y) w(y,t)dy  + \alpha \int_0^x k^a(x,y) v(y,t)dy   \nonumber \\
     & =  \tilde{w}_{xx}(x,t) + \int_0^x k^a_{xx}(x,y)w(y,t)dy +k^a_x(x,x) w(x,t) \nonumber \\
     & + \frac{d}{dx} k^a(x,x) w(x,t)+ k^a(x,x) w_x(x,t) -\rho w(x,t) \nonumber \\
     &+ \alpha v(x,t) . \label{new-eq}
\end{align}
Since $k^a_y(x,0)=0$, then adding and subtracting $(c_2+\rho) \tilde{w}(x,t) $ to the right-hand-side of \eqref{new-eq}
\begin{align*}
     & \tilde{w}_{t}(x,t)   
         =  \tilde{w}_{xx}(x,t)  - (c_2+\rho) \tilde{w}(x,t)+ \alpha v(x,t)    \nonumber \\
         &- \alpha \int_0^x k^a(x,y) v(y,t)dy + (2\frac{d}{dx} k^a(x,x) + c_2) w(x,t) \nonumber \\
         & + \int_0^x [k^a_{xx}(x,y) - k^a_{yy}(x,y)  - c_2  k^a(x,y) ] w(y,t)dy   .
\end{align*}
Since $k^a(x,y)$ is given by \eqref{26}, the previous equation reduces to \eqref{19}. Also,
\begin{align*}
    \tilde{w}_x(0,t)=w_x (0,t) -k^a(0,0) w(0,t) = 0,
\end{align*}
 and  the other boundary condition on $w(x,t)$  holds by using \eqref{27}. Equation \eqref{20} can be obtained by referring to \eqref{inverse of transformation 1}.
 $\hfill \square$
 
Next, we provide conditions that ensure the exponential stability of the target system. First, we need  the following lemma,  which provides bounds on the induced $L^2$-norms of the kernel functions $k^a(x,y)$ and $\ell^a(x,y)$ .

\begin{lem} \label{bounds on k, l}
    The $L^2$-norms of  $k^a(x,y)$ and $\ell^a(x,y)$  are bounded by
    \begin{align}
         \| k^a\| \leq & \sqrt{\frac{c_2\pi}{8}} \; \left( erfi (\sqrt{\frac{c_2}{2}} )  erf (\sqrt{\frac{c_2}{2}} ) \right)^{\frac{1}{2}}, \label{29} \\
        \| \ell^a\| \leq & \sqrt{\frac{c_2\pi}{8}} \; \left( erfi (\sqrt{\frac{c_2}{2}} )  erf (\sqrt{\frac{c_2}{2}} ) \right)^{\frac{1}{2}}, \label{30} 
    \end{align}
   where $erfi (x) = \frac{2}{\sqrt{\pi}} \int_0^x e^{\xi^2} d \xi$,   $ \;erf (x) = \frac{2}{\sqrt{\pi}} \int_0^x e^{-\xi^2} d \xi$ . 
\end{lem}
 \noindent \textbf{Proof.} To prove relation \eqref{29}, we recall the expression for the kernel $k^a(x,y)$ given in \eqref{k-kernel}.
 We set $z=\sqrt{c_2(x^2-y^2)}$, then
 \begin{align*}
     k^a(x,y) =& \frac{-c_2}{z}  x \sum_{m=0}^{\infty} \left( \frac{z}{2}\right)^{2m+1} \frac{1}{m! m+1!} \\
    % =&\frac{-c_2}{z}  \; x \; \frac{z}{2}  \sum_{m=0}^{\infty} \left( \frac{z}{2}\right)^{2m} \frac{1}{m! m+1!}\\
     =& \frac{-c_2}{2} x \sum_{m=0}^{\infty} \frac{(z^2/4)^m}{m!} \frac{1}{ m+1!} \\
     \leq &  \frac{-c_2}{2} x \sum_{m=0}^{\infty} \frac{(z^2/4)^m}{m!}.
 \end{align*}
 Thus the induced $L_2$-  norm of $k^a(x,y)$ is bounded by 
 \begin{align*}
     & \| k^a(x,y) \| \leq \frac{c_2}{2} \|x\| \|e^{\frac{z^2}{4}} \| \\
     & \leq  \frac{c_2}{2} \|x\| \|e^{\frac{c_2x^2}{4}} \|  \|e^{\frac{-c_2y^2}{4}}\|  \\
     & \leq   \sqrt{\frac{c_2\pi}{8}} \; \left( erfi (\sqrt{\frac{c_2}{2}} )  erf (\sqrt{\frac{c_2}{2}} ) \right)^{\frac{1}{2}}.
 \end{align*}
 Similarly, one can prove \eqref{30} by referring back to \eqref{28}. 
 \begin{align*}
     \ell^a(x,y) &=  \frac{-c_2}{z} x  \sum_{m=0}^{\infty}  (-1)^m\left( \frac{z}{2}\right)^{2m+1} \frac{1}{m! m+1!} \\ 
     &\leq   \frac{c_2}{z} x  \sum_{m=0}^{\infty}  \left( \frac{z}{2}\right)^{2m+1} \frac{1}{m! m+1!} ,
\end{align*}
and the $L_2$-norm of $l(x,y)$ is bounded by 
\begin{align*}
     &\|\ell^a(x,y)\|\leq  \sqrt{\frac{c_2\pi}{8}} \; \left( erfi (\sqrt{\frac{c_2}{2}} )  erf (\sqrt{\frac{c_2}{2}} ) \right)^{\frac{1}{2}}. 
 \end{align*}
 \hfill $\square$ 
 
 The following lemma will be needed to show  stability of the target system.

\begin{lem} \label{lemma} Let $\gamma>0$. The states  of the target system (\ref{19})-(\ref{22}) satisfy 
    \begin{eqnarray}
    \| v(x,t)\| \leq \frac{|\beta|}{\gamma} (1+\|\ell^a\|) \| \tilde{w}\|. \label{tt9}
\end{eqnarray}
\end{lem}
\noindent \textbf{Proof.} Multiply equation \eqref{20} with $v(x,t)$ and integrate from $0$ to $1$, 
\begin{align*}
    & 0 = \int_0^1 v_{xx}(x,t) v(x,t) dx - \gamma \int_0^1 v^2(x,t)dx  \nonumber \\ 
         &  + \beta \int_0^1 \tilde{w}(x,t) v(x,t) dx   + \beta  \int_0^1 v(x,t) \int_0^x \ell^a(x,y) \nonumber \\ 
         & \times \tilde{w}(y,t) dy dx.
\end{align*}
Thus
\begin{align}
     &\gamma \int_0^1 v^2(x,t)dx  \leq  \beta \int_0^1 \tilde{w}(x,t) v(x,t) dx  \nonumber \\
     &+ \beta  \int_0^1 v(x,t) \int_0^x \ell^a(x,y) \tilde{w}(y,t) dy dx. \label{tt66}
\end{align}
Bounding the  terms on the right-hand side of inequality \ref{tt66} using Cauchy-Schwartz leads to \eqref{tt9}.
 \hfill $\square$ \\
\smallskip
\begin{thm}\label{thm-stab}
    The target system $\eqref{19}-\eqref{22}$ is exponentially stable if   
\begin{align}
   c_2  + \rho > & \frac{|\alpha \beta|}{\gamma} (1+\|\ell^a\|)  (1+\|k^a\|). \label{c_1} 
\end{align} 
\end{thm}
\noindent \textbf{Proof.}  Define the Lyapunov function candidate,
\begin{align*}
    V(t) =& \frac{1}{2} \int_0^1 \tilde{w}^2(x,t)dx = \frac{1}{2} \| \tilde{w}(x,t) \|^2  .
\end{align*}
Taking the time derivative of $V(t)$,
\begin{align}
    &\dot{V}(t) = \int_0^1 \tilde{w}(x,t) \tilde{w}_t(x,t) dx  \nonumber \\
  % &= \int_0^1  \tilde{w}(x,t) [ \tilde{w}_{xx}(x,t) - (c_2+\rho) \tilde{w}(x,t)+ \alpha v(x,t) \nonumber \\
  % &- \alpha \int_0^x k^a(x,y) v(y,t) dy ]dx \nonumber \\
    &\leq  -(c_2+\rho) \int_0^1 \tilde{w}^2(x,t) dx + \alpha \int_0^1 \tilde{w}(x,t) v(x,t) dx\nonumber \\
    &- \alpha \int_0^1 \tilde{w}(x,t)\int_0^x k^a(x,y) v(y,t) dydx . \label{31}
\end{align}
Using Cauchy-Schwartz inequality, we estimate the term of the  right-hand-side of inequality \eqref{31} as follows.
\begin{align}
    \alpha \int_0^1 \tilde{w}(x,t) v(x,t) dx \leq & |\alpha| \|\tilde{w}\| \|v\| \nonumber \\
      \leq & \frac{|\alpha||\beta|}{\gamma} (1+\|\ell^a\|) \| \tilde{w}\|^2 , \label{32}
\end{align}
and
\begin{align}
  & - \alpha \int_0^1 \tilde{w}(x,t) \int_0^x k^a(x,y) v(y,t) dy dx \nonumber \\
  %&\leq  |\alpha|  \int_0^1 |\tilde{w}(x,t)|\int_0^1 |k^a(x,y)| |v(y,t)| dy dx\nonumber \\
    &\leq  |\alpha| \|k^a\| \| \tilde{w}\| \| v \| \nonumber \\
     &\leq  \frac{|\alpha \beta|}{\gamma}   \|k^a\| (1+\|\ell^a\|)    \| \tilde{w}\|^2 . \label{33}
\end{align}
Subbing \eqref{32} and \eqref{33} in \eqref{31}, 
\begin{align}
   \dot{V}(t) \leq&   -\left((c_2+\rho)  - \frac{|\alpha \beta|}{\gamma} (1+\|\ell^a\|) (1+\|k^a\|)\right) \| \tilde{w}\|^2.\label{34}
\end{align}
Setting
\begin{align}
c_3=& (c_2+\rho)  - \frac{|\alpha \beta|}{\gamma} (1+\|\ell^a\|) (1+\|k^a\|), \label{bound on the decay rate}
\end{align}
then inequality \eqref{34} implies that  $   V(t) \leq  e^{-2 c_3  t} V(0).$ If the parameter $c_2$ is chosen such that \eqref{c_1} is satisfied, then $V(t)$ decays exponentially as $t \rightarrow \infty$, and so does $\| \tilde{w}(x,t)\|$. By means of lemma \eqref{lemma}, the state $v(x,t)$ is asymptotically stable.  Recalling  that the operator $(\partial_{xx} - \gamma I)$ is boundedly invertible, then the elliptic equation \eqref{20} implies that
\begin{align*}
   &v(x,t) 
   \nonumber \\
   & =   (  \gamma I-\partial_{xx})^{-1} \left( \beta \tilde{w}(x,t)   + \beta  \int_0^x \ell^a (x,y) \tilde{w}(y,t)  dy \right).
\end{align*}
Substituting for $v(x,t)$ in the parabolic equation \eqref{19} leads to a system described by the $\Tilde{w}(x,t)$ only. Hence, the exponential stability of the coupled system follows from the exponential stability of the state $\tilde{w}(x,t)$ . 
\hfill   $\square$

The decay rate of the target system is bounded by \eqref{bound on the decay rate}. The following theorem  is now immediate.

\begin{thm}
System \eqref{1}-\eqref{4} is exponentially stable if the control signal is  
\begin{align}
        u(t)=&\int_0^1 k^a_x(1,y) w(y,t)dy+k^a(1,1) w(1,t), \label{control signal stab}  
\end{align}
with $k^a(x,y)$  as in lemma\ref{Wellposedness}, and parameter $c_2$ satisfies 
\begin{align}
   &c_2 + \rho \nonumber \\
   &>  \frac{|\alpha \beta|}{\gamma}  \left[ 1+\sqrt{\frac{c_2\pi}{8}} \; ( erfi (\sqrt{\frac{c_2}{2}} ) )^{\frac{1}{2}} ( erf (\sqrt{\frac{c_2}{2}} ) )^{\frac{1}{2}}  \right]^2.  \label{bound on c1} 
\end{align} 
\end{thm}
\noindent \textbf{Proof.}
Since $c_2 $ is given by \eqref{bound on c1}, it follows from \eqref{thm-stab} and \eqref{bounds on k, l} that the target system \eqref{19}-\eqref{22} is exponentially stable. It follows from Theorem \ref{thm1} that with $u(t)$ given as in \eqref{control signal stab}, there is an invertible state transformation between system \eqref{1}-\eqref{4} and the exponentially stable target system  \eqref{19}-\eqref{22}. The conclusion is now immediate. $\hfill \square$

\begin{figure}[H]
\begin{center}
\textbf{Illustration of the restriction \eqref{bound on c1}  on $c_2$  }\par\medskip
{\includegraphics[scale=0.3]{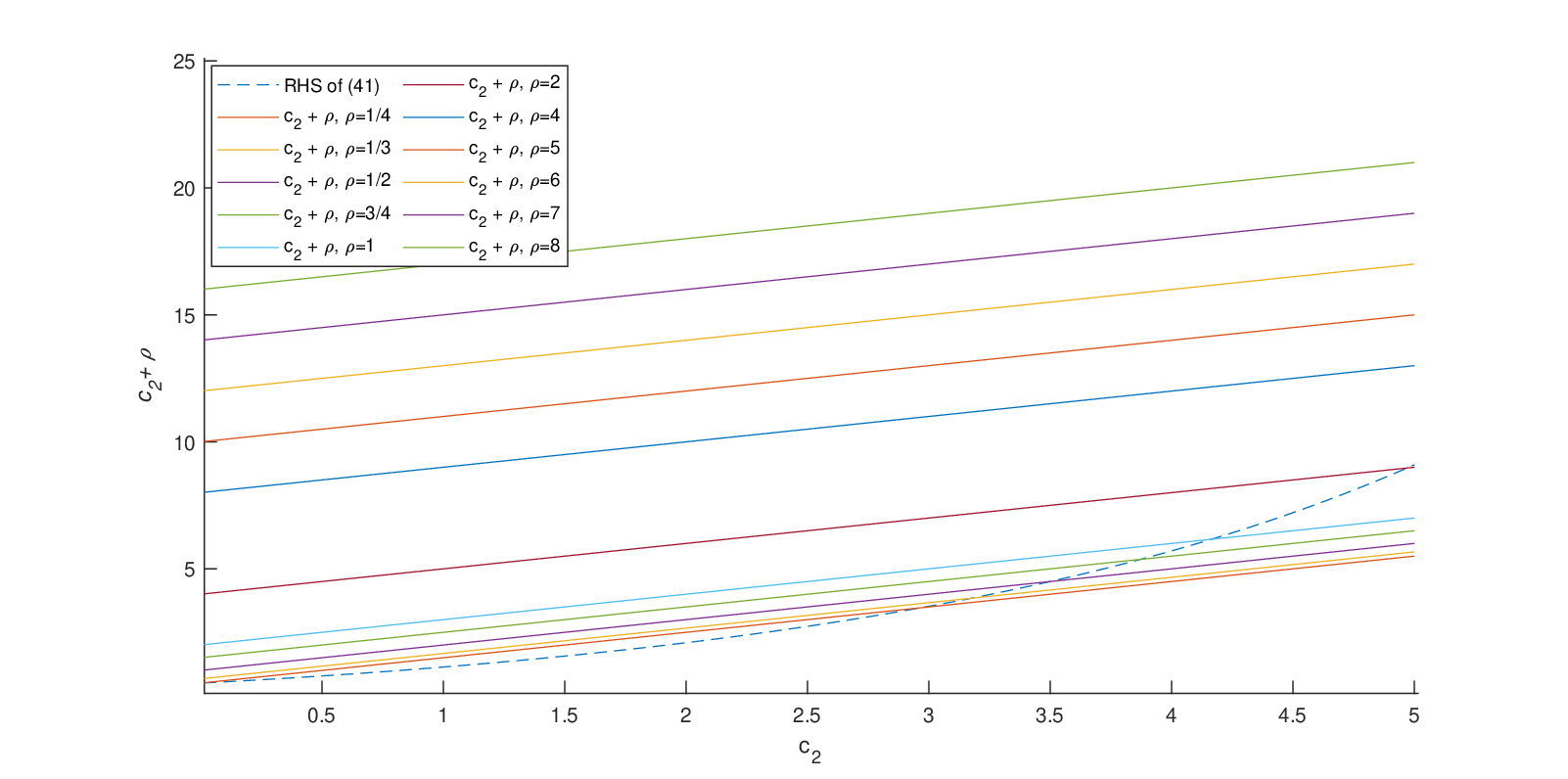} }
\end{center}
\caption{A comparison between the right-hand-side of \eqref{bound on c1} as a function of $c_2$ against several straight lines  $c_2+\rho $ for different values of $\rho$, where the other parameters are fixed as $\beta=\gamma=1$, $\alpha=0.5 .$ The  right-hand-side of \eqref{bound on c1} is described using a dashed line(- - -). Target system \eqref{19}-\eqref{22} is exponentially stable for  values of $c_2$ at which  the straight line $c_2+\rho$, for some $\rho$,  is above the dashed line(- - -). The figure showcases the restrictive nature associated with  condition \eqref{bound on c1}. The parameter $\rho$ has to be  large %for the criteria outlined in \eqref{bound on c1} to be fulfilled 
or  the coupling factor $\alpha \beta$ has to be small allowing for the inequality \eqref{bound on c1} to be fulfilled.   }
\label{c2+rho-comparison}
\end{figure}

Figure\ref{c2+rho-comparison} illustrates the  restrictiveness of the  criterion \eqref{bound on c1}. This figure gives a comparison between the right-hand-side of inequality \eqref{bound on c1} and different straight lines $c_2+\rho$ for various values of $\rho$ while setting  $\gamma=\beta=1$ and $\alpha=0.5 .$ The dashed line describes the right-hand-side of inequality \eqref{bound on c1}, whereas the straight lines present straight lines $c_2+\rho$, for different values of $\rho$.  For some $\rho$, if  values of $c_2$ are such that the dashed line (- - -) is below the straight line $c_2 + \rho$, bound \eqref{bound on c1}  is  fulfilled, and hence stability of the target system $\eqref{19}-\eqref{22}$ follows. % Alternatively, the bound is satisfied  if $\rho$ is sufficiently large compared to $\frac{\alpha \beta}{\gamma}$  implying, in this case, that the system is already exponentially stable by looking at inequality \eqref{16}. 

\subsection{Numerical simulations} \label{numerical simulatino1}
%%%%%%%%%%%%%%%%%%%%%%%%%%%%%%%%%%%%%%%%%%%%%%%%%%%% Numerical simulations-stab

\begin{figure}[H]
\begin{center}
 \textbf{Comparison of  open-loop and closed loop dynamics }\par\medskip
\subfloat[uncontrolled $w(x,t)$]{\includegraphics[scale=0.25]{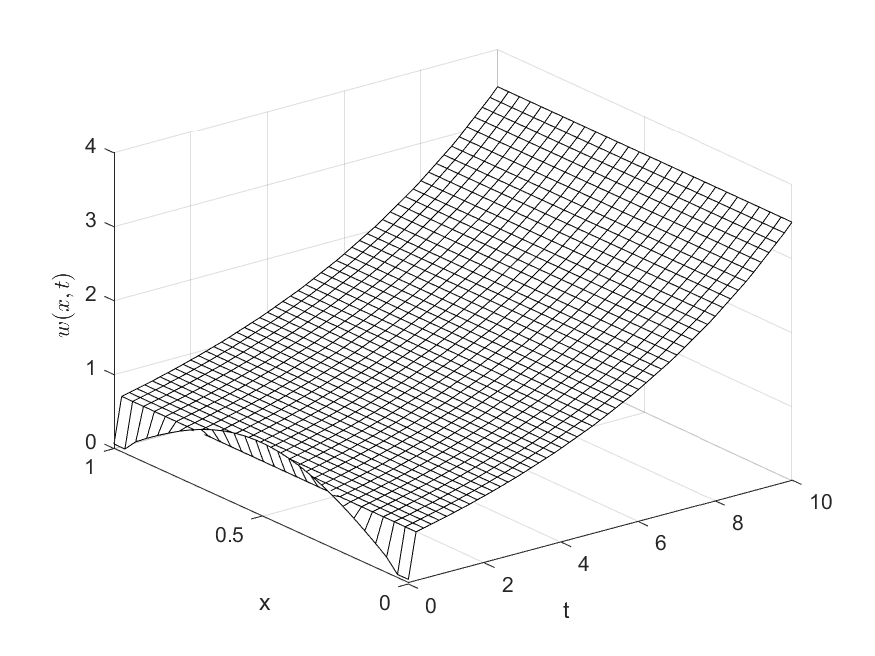}}
\hfill
\subfloat[uncontrolled $v(x,t)$]{\includegraphics[scale=0.25]{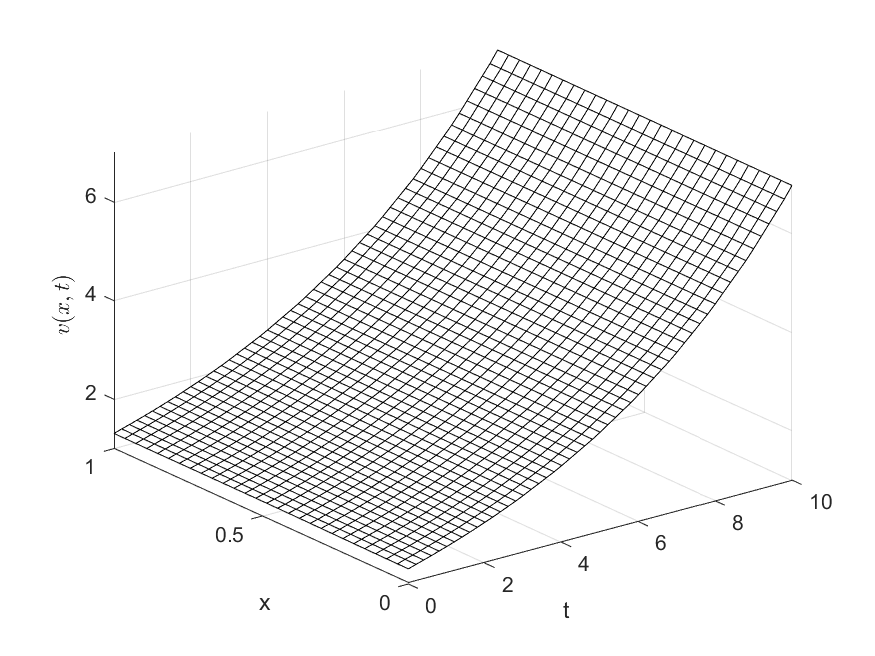}} \\
\subfloat[controlled $w(x,t)$]{\includegraphics[scale=0.25]{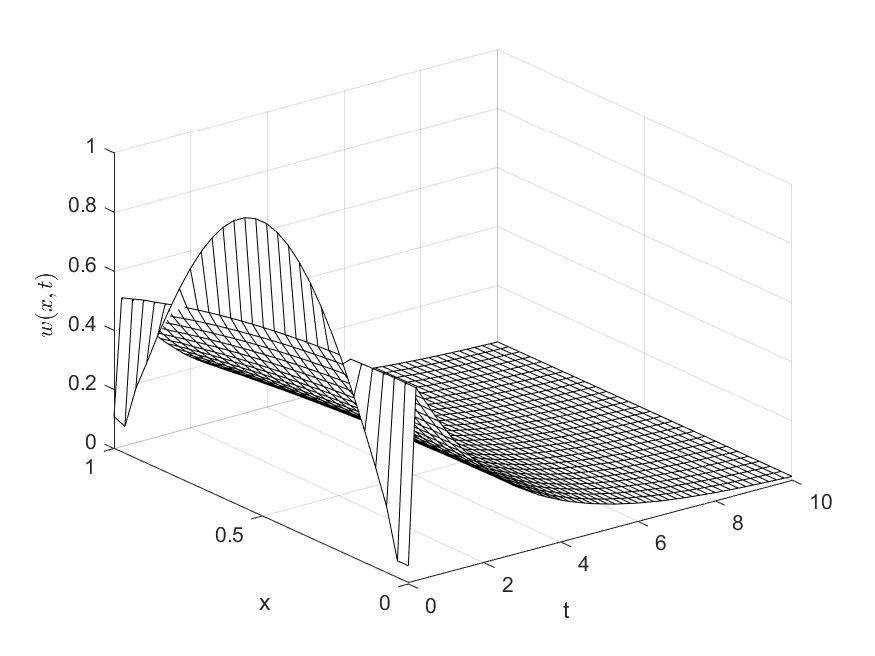}}
\hfill
\subfloat[controlled $v(x,t)$]{\includegraphics[scale=0.25]{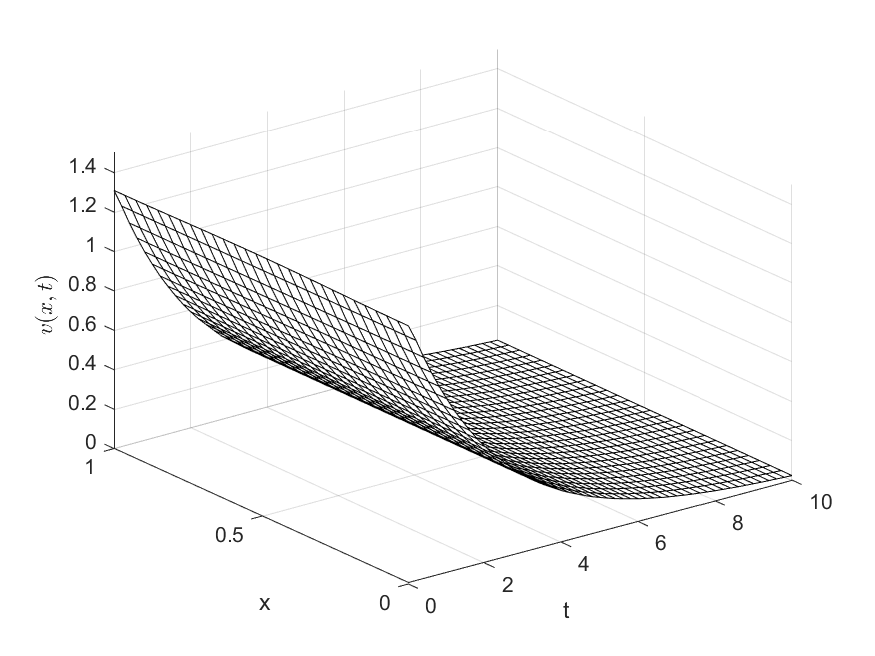}}
\end{center}
\caption{A 3D landscape of the dynamics of a coupled parabolic-elliptic system (\ref{1})-(\ref{4}) with initial condition  $w_0=\sin(\pi x)$, without and with  control.  The parameters of the system are  $\gamma=\frac{1}{4}$, $\rho=\frac{1}{3}$,  $\alpha =\frac{1}{4} $ , $\beta=\frac{1}{2}.$  The  uncontrolled system is unstable with this choice of parameters. The backstepping state-feedback control gain is  $c_2= 1.2-\rho$ which meets  the stability  condition   \eqref{bound on c1}. The control causes the solutions of the system to decay to the steady-state solution as $t \to \infty$. }
\label{fig1}
\end{figure}

\begin{figure}[H]
\begin{center}
 \textbf{Comparison of open-loop and closed-loop  $L_2$-norms }\par\medskip
\subfloat[$L_2$-norm $w(x,t)$]{\includegraphics[scale=0.25]{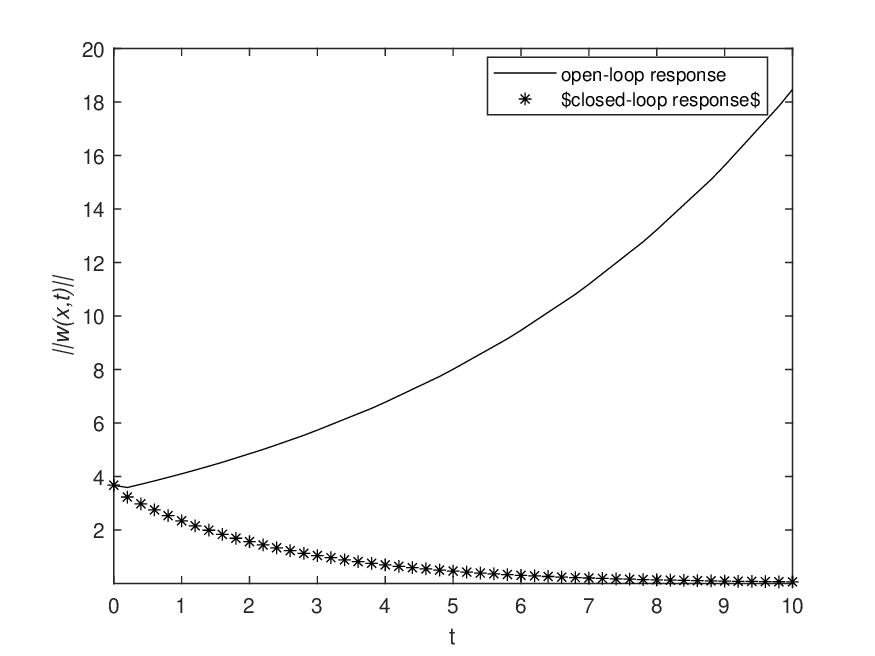}}
\hfill
\subfloat[$L_2$-norm  $v(x,t)$]{\includegraphics[scale=0.25]{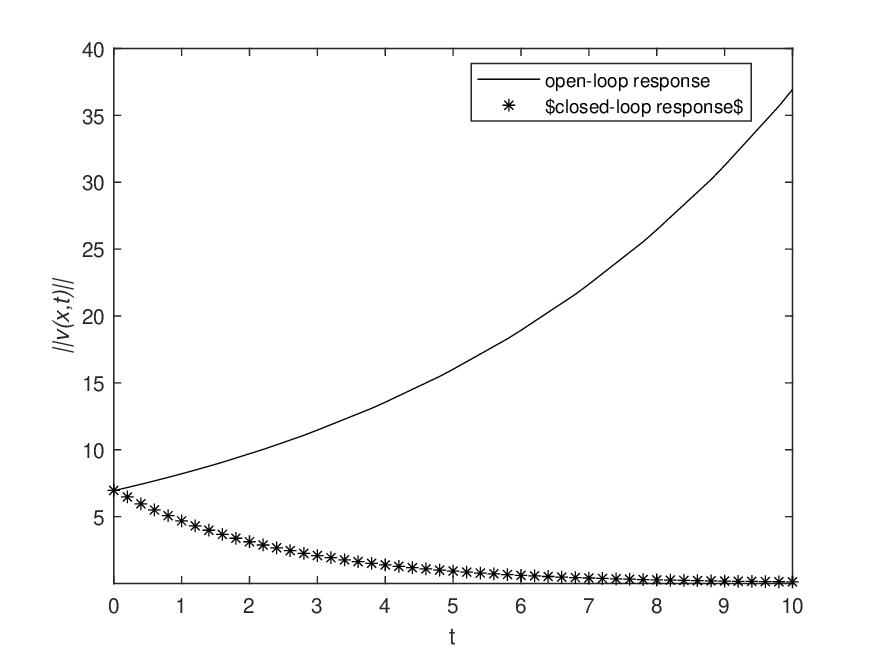} }
\end{center}
\caption{ A comparison between the $L_2$-norm of the solutions $w(x,t)$ and $v(x,t)$ for the uncontrolled and controlled systems with initial condition  $w_0=\sin(\pi x).$ The system parameters  are $\gamma=\frac{1}{4}$, $\rho=\frac{1}{3}$,  $\alpha =\frac{1}{4} $, $\beta=\frac{1}{2}$ so the uncontrolled system is unstable. In  the absence of control, the $L_2$-norm of the states grows.  The backstepping state-feedback control gain is given by $ c_2  = 1.2- \rho ;$  stability condition \eqref{bound on c1} is satisfied.  The $L_2$-norm of the solution of the controlled system  decays to zero. }
\label{fig3}
\end{figure}

The solutions of system \eqref{1}-\eqref{4}, both controlled and uncontrolled, were simulated numerically  using a finite-element approximation in COMSOL Multiphysic software. The finite-element method (FEM)  with linear splines was used to approximate the coupled equations by a system of DAEs. The spatial interval was  divided into $27$ subintervals. Also, time was discretized by  a time-stepping algorithm called  generalized alpha with time-step= 0.2. We set $\gamma=\frac{1}{4}$, $\rho=\frac{1}{3}$,  $\alpha =\frac{1}{4} $ and $\beta=\frac{1}{2}.$  For these parameter values, the system is unstable. Figure\ref{fig1} presents the dynamics of the  states $w(x,t)$ and $v(x,t)$ in the absence of the control with initial condition    $w(0)=\sin (\pi x ) .$  The system was first controlled with the controller resulting from the choice of parameter $c_2=1.2 -\rho $ which satisfies inequality  \eqref{bound on c1} and thus stability of the controlled system is guaranteed. This is illustrated in Figure\ref{fig1}. As predicted by the theory, the dynamics of the system decay to zero with time. A comparison between the $L_2$-norm of  both states $w(x,t)$ and $v(x,t)$ before and after applying the control input is given in Figure \ref{fig3}.

\begin{comment}
\begin{figure}[H]
\centering
  \textbf{Dynamics of closed-loop \eqref{1}-\eqref{4}}\par\medskip
\begin{center}
\subfloat[$w(x,t)$]{\includegraphics[scale=0.25]{w-after control alpha=0.25....eps}}
\hfill
\subfloat[$v(x,t)$]{\includegraphics[scale=0.25]{v-after control alpha=0.25....eps}}
\end{center}
\caption{A 3D landscape of the dynamics of a  coupled parabolic-elliptic system (\ref{1})-(\ref{4}) after applying the control \eqref{control signal stab}. The dynamics was simulated by employing finite-element approximation using a linear basis. The system parameters were set as  $\gamma=\frac{1}{4}$, $\rho=\frac{1}{3}$,  $\alpha =\frac{1}{4} $ , $\beta=\frac{1}{2}$ and the initial condition  $w_0=\sin(\pi x)$. The control gain is given by $c_2= 1.2-\rho$. Consequently, the system meets  the stability  condition   \eqref{bound on c1}. A comparison between these simulations and Figure\ref{fig1} clearly illustrates the impact of  control \eqref{control signal stab}. The control forces the solutions of the system to decay to the steady-state solution as $t \to \infty$.}
\label{fig2}
\end{figure}
\end{comment}

\section{Observer design} \label{section3}

In the previous section, the control input was designed based on the assumption that the state of system \eqref{1}-\eqref{4} is known. The focus of this section is to design an observer that  estimates the state of the  parabolic-elliptic system. 

\begin{comment}
\chgk{Really? The designs look similar...Such design  relies on two critical factors: (i) the type of the  measurements provided, whether they  are given for the  parabolic state , the elliptic state or  both, and (ii) the specific location designated for the given measurements. }

The construction of the observer depends on these factors, as they have direct impact on the analysis taken to establish the exponential stability of the observation error dynamics. Although the latter is demonstrated through the exploitation of a backstepping approach for all the different cases, the transformation employed in the backstepping approach will be tailored depending on the aforementioned factors accordingly.  In turn each customized transformation enables the selection of appropriate output injections that helps in achieving the observation error dynamics' stabilization.
\end{comment}
We first design an observer using  two measurements   $v(1, t)$, and $w(1, t)$. As for the situation when two controls can be used, the design is fairly straightforward. The situation becomes more intricate with only one measurement is available. Then, as for the single control situation covered in the previous section, a bound must be satisfied.  

\subsection{Both $w(1,t)$ and $v(1,t)$ are available}\label{subsection1}
The objective is to design an observer when the available  measurements of system \eqref{1}-\eqref{4} are $w(1,t)$ and $v(1,t)$. We propose the following observer for system \eqref{1}-\eqref{4}
\begin{subequations}  \label{equation56}
\begin{align}
     \hat{w}_t(x,t) =&  \hat{w}_{xx}(x,t) - \rho \hat{w}(x,t)  + \alpha \hat{v}(x,t)   \nonumber \\
     &+ \eta_1(x) [ w(1,t) - \hat{w}(1,t)],  \label{equation56a}\\
    0 =&  \hat{v}_{xx}(x,t)  - \gamma \hat{v}(x,t) + \beta \hat{w}(x,t)    \nonumber \\
     &  + \eta_2(x) [ v(1,t) - \hat{v}(1,t)],  \label{equation56b}\\
   \hat{w}_x(0,t) =& 0 ,  \quad \hat{w}_x(1,t) =u(t)+ \eta_3 [ w(1,t)- \hat{w}(1,t)],&  \label{equation56c} \\
   \hat{v}_x (0,t) =&0,    \quad \hat{v}_x(1,t) =\eta_4 [ v(1,t)- \hat{v}(1,t)],  &
 \label{equation56d}
\end{align}
\end{subequations}
Two in-domain output injection functions $\eta_1(x)$ and $\eta_2(x)$, and two boundary injections  values $\eta_3$ and $\eta_4$ are to be designed. We will need the following lemma. 
\begin{lem} \cite{smyshlyaev2005backstepping} \label{Wellposednesswv1}
    The hyperbolic partial differential equation
 \begin{subequations} \label{equation58}
\begin{align}
    & k^b_{xx}(x,y) - k^b_{yy}(x,y)  + o_2 k^b(x,y) = 0, \quad 0<x<y<1 \label{equation58a}  \\
    & k^b_x(0,y)= 0, \quad   k^b(x,x)= -\frac{1}{2} o_2 x,   \label{equation58b}
\end{align}
\end{subequations}
 has a continuous unique solution. Here $o_2=o_1-\rho$, and $o_2>0$.
\end{lem}
Define the error states
\begin{subequations}  \label{error states}
    \begin{align}
   &e^w(x,t) =  w(x,t) -  \hat{w}(x,t), \label{error states1} \\
    &e^v(x,t) =  v(x,t) -  \hat{v}(x,t).     \label{error states2}
  \end{align} 
\end{subequations}
The observer error dynamics satisfy 
\begin{subequations} \label{equation57}
    \begin{align}
        e^w_t(x,t) =& e^w_{xx}(x,t) - \rho e^w(x,t) + \alpha e^v(x,t) \nonumber \\
        & - \eta_1(x) e^w(1,t), \label{equation57a} \\
        0=& e^v_{xx} (x,t) - \gamma e^v (x,t) + \beta e^w(x,t)   \nonumber \\
        & - \eta_2(x) e^v(1,t),\label{equation57b} \\
        e^w_x(0,t)=& 0,  \quad e^w_x(1,t)=- \eta_3 e^w(1,t),& \label{equation57c} \\
         e^v_x(0,t)=& 0,   \quad e^v_x(1,t)=- \eta_4 e^v(1,t), &\label{equation57d}
    \end{align}
\end{subequations}

A backstepping approach is used to  select $\eta_1(x), \; \eta_2(x),$ $ \eta_3, \; \eta_4$  so that the error system \eqref{equation57d} is exponentially stable.
We introduce the  target system
\begin{subequations} \label{equation47}
    \begin{align}
        e^{\tilde{w}}_t(x,t) =&  e^{\tilde{w}}_{xx}(x,t) - (o_2+\rho) e^{\tilde{w}}(x,t) + \alpha e^{\tilde{v}}(x,t), \label{equation47a}  \\
        0=& e^{\tilde{v}}_{xx} (x,t) - (o_2 +\gamma)  e^{\tilde{v}} (x,t) + \beta  e^{\tilde{w}}(x,t) , \label{equation47b}\\
         e^{\tilde{w}}_x(0,t)=& 0, \qquad  e^{\tilde{w}}_x(1,t)=0, \label{equation47c}\\
         e^{\tilde{v}}_x(0,t)=& 0, \qquad e^{\tilde{v}}_x(1,t)=0. \label{equation47d}
  \end{align} 
\end{subequations}
A pair of state transformations
\begin{subequations} \label{transformationwv1}
\begin{align} 
        e^w(x,t)=&e^{\tilde{w}}(x,t)-\int_x^1 k_1(x,y) e^{\tilde{w}}(y,t)dy, \label{equation59} \\
        e^v(x,t)=&e^{\tilde{v}}(x,t)-\int_x^1 k_2(x,y) e^{\tilde{v}}(y,t)dy, \label{equation60}
\end{align}
\end{subequations}
that transform the target system \eqref{equation47} into \eqref{equation57} are needed.  

\begin{thm} \label{Theorem 11wv1} If   $k_1(x,y)=k_2(x,y)=k^b(x,y)$ where $k^b(x,y)$ satisfies \eqref{equation58}, and if the output injections  are 
\begin{align}
&\eta_1(x) = \eta_2(x)= -k^b_y(x,1), \label{outputinjection2.1}  \\
&\eta_3 =  \eta_4=-k^b(1,1),\label{outputinjection2.2}  
\end{align}
then transformations \eqref{equation59} and \eqref{equation60}  convert the target system \eqref{equation47} into the original  error dynamics  \eqref{equation57}.
\end{thm}   
\noindent \textbf{Proof.} We first take the spatial   derivatives of \eqref{equation59}.
\begin{align}
& e_x^w(x,t)=e^{\tilde{w}}_{x}(x,t)-\int_x^1 k_{1x}(x,y) e^{\tilde{w}}(y,t) dy \nonumber \\
& + k_1(x,x) e^{\tilde{w}}(x,t) \label{one derivative wv1}\\
          &e^w_{xx}(x,t)= e^{\tilde{w}}_{xx}(x,t)-\int_x^1 k_{1xx}(x,y) e^{\tilde{w}}(y,t)dy  \nonumber \\
& + k_{1x}(x,x) e^{\tilde{w}}(x,t) + \frac{d}{dx} k_1(x,x) e^{\tilde{w}}(x,t) \nonumber \\ 
     & + k_1(x,x) e^{\tilde{w}}_x(x,t). \hfill \label{equation64}
\end{align}
     Taking the time   derivative of \eqref{equation59} and integrating by parts
\begin{align}
         &e^w_{t}(x,t)= e^{\tilde{w}}_{t}(x,t)-\int_x^1 k_1(x,y) e^{\tilde{w}}_t(y,t)dy \nonumber \\
    % &= e^{\tilde{w}}_t(x,t) - \int_x^1 k_1(x,y) [ e^{\tilde{w}}_{yy}(y,t) - (o_2+\rho) e^{\tilde{w}}(y,t)  \nonumber \\ + \alpha e^{\tilde{v}}(y,t)]dy \nonumber \\
      &=   e^{\tilde{w}}_t(x,t)   + (o_2+\rho) \int_x^1 k_1(x,y) e^{\tilde{w}}(y,t)dy   - k_1(x,1) \nonumber \\
      & \times e^{\tilde{w}}_x(1,t)  - \alpha \int_x^1 k_1(x,y) e^{\tilde{v}}(y,t)dy   +   k_1(x,x) e^{\tilde{w}}_x(x,t)     \nonumber \\ 
     & + k_{1y}(x,1) e^{\tilde{w}}(1,t) -k_{1y}(x,x) e^{\tilde{w}}(x,t) \nonumber \\
     &-   \int_x^1 k_{1yy} (x,y)  e^{\tilde{w}}(y,t) dy.
         \label{equation65}
\end{align}
We rewrite the right-hand-side  of the parabolic equation \eqref{equation57a} of the error dynamics as
\begin{align}
      &e^w_t(x,t) - e^w_{xx}(x,t) + \rho e^w(x,t) - \alpha e^v(x,t) \nonumber \\ 
      &+ \eta_1(x) e^w(1,t)=0. \label{parabolic-rearrange}
\end{align}
Substituting \eqref{equation64} and \eqref{equation65} in \eqref{parabolic-rearrange}, then  the left-hand-side of \eqref{parabolic-rearrange} is 
 \begin{align}
      &  \text{(L.H.S)}_1= e^{\tilde{w}}_t(x,t) -e^{\tilde{w}}_{xx}(x,t) - \int_x^1 k_{1yy}(x,y)  \nonumber \\
      & \times e^{\tilde{w}}(y,t) dy  + (o_2+\rho)\int_x^1 k_1(x,y) e^{\tilde{w}}(y,t)dy \nonumber \\
      & +\int_x^1 k_{1xx}(x,y)e^{\tilde{w}}(y,t)dy  - k_{1y}(x,x) e^{\tilde{w}}(x,t) \nonumber \\
      & - k_{1x}(x,x) e^{\tilde{w}}(x,t)  - \frac{d}{dx} k_1(x,x) e^{\tilde{w}}(x,t) - k_1(x,x) \nonumber \\
      & \times e^{\tilde{w}}_x(x,t) + k_1(x,x) e^{\tilde{w}}_x(x,t) + \rho e^{\tilde{w}}(x,t) - \alpha e^{\tilde{v}}(x,t)  \nonumber \\ 
     & - \rho \int_x^1 k_1(x,y) e^{\tilde{w}}(y)dy  + \eta_1(x) e^{\tilde{w}}(1,t) + k_{1y}(x,1) \nonumber \\ 
     &\times e^{\tilde{w}}(1,t) -   k_1(x,1) e^{\tilde{w}}_x(1,t) -\alpha \int_x^1 k_1(x,y)  e^{\tilde{v}}(y,t)dy  \nonumber \\
     &  + \alpha \int_x^1 k_2(x,y) e^{\tilde{v}}(y,t)dy , \label{rearranging termswv1} \\
     &  \text{(R.H.S)}_1=0. \label{right-handside2}
 \end{align}
 Adding and subtracting the term $(o_2+\rho) e^{\tilde{w}}(x,t)$ to the right-hand-side of \eqref{rearranging termswv1}  
\begin{align}
 &  \text{(L.H.S)}_1= e^{\tilde{w}}_t(x,t)   - e^{\tilde{w}}_{xx}(x,t)  +o_2  e^{\tilde{w}}(x,t)  - \alpha  e^{\tilde{v}}(x,t)  \nonumber \\
 & -   \int_x^1 [-o_2 k_1(x,y)  - k_{1xx}(x,y)+  k_{1yy}(x,y)  ] e^{\tilde{w}}(y,t) dy \nonumber \\
& -   k_1(x,1) e^{\tilde{w}}_x(1,t) - (2  \frac{d}{dx} k(x,x)  + o_2)    e^{\tilde{w}}(x,t)  \nonumber \\
& + (\eta_1(x) + k_{1y}(x,1)) e^{\tilde{w}}(1,t)  -\alpha \int_x^1 k_1(x,y)  e^{\tilde{v}}(y,t)dy \nonumber \\
 & + \alpha \int_x^1 k_2(x,y) e^{\tilde{v}}(y,t)dy .\label{equationwv1}
\end{align}
Using the  boundary condition $ e^{\tilde{w}}_x(1,t)=0$,  if $k_1(x,y)=k_2(x,y)=k^b(x,y)$    then equation \eqref{equationwv1} reduces to 
\begin{align}
 & \text{(L.H.S)}_1=  e^{\tilde{w}}_t(x,t)   - e^{\tilde{w}}_{xx}(x,t)  +(o_2+\gamma)  e^{\tilde{w}}(x,t) \nonumber \\
 & - \alpha  e^{\tilde{v}}(x,t)  -   \int_x^1 [-o_2 k^b(x,y)  - k^b_{xx}(x,y)+  k^b_{yy}(x,y)  ] \nonumber \\ & \times e^{\tilde{w}}(y,t) dy  - (2  \frac{d}{dx} k^b(x,x)  + o_2)    e^{\tilde{w}}(x,t) \nonumber\\
 & + (\eta_1(x) + k^b_y(x,1)) e^{\tilde{w}}(1,t)  . \label{new-equation-2}
\end{align}
If $k^b(x,y)$ satisfies \eqref{equation58} and   $\eta_1(x) =k^b_y(x,1)$, then \eqref{new-equation-2} becomes
\begin{align*}
 & \text{(L.H.S)}_1=  e^{\tilde{w}}_t(x,t)   - e^{\tilde{w}}_{xx}(x,t)  +(o_2+\gamma)  e^{\tilde{w}}(x,t)  \nonumber \\
 &- \alpha  e^{\tilde{v}}(x,t) .
\end{align*}
Referring to \eqref{right-handside2} and \eqref{equation47a},
\begin{align*}
    \text{(L.H.S)}_1=0=\text{(R.H.S)}_1.
\end{align*}
Hence the state transformation \eqref{transformationwv1} transforms the parabolic equation \eqref{equation47a} into \eqref{equation57a}. Referring to \eqref{one derivative wv1}, we apply transformation \eqref{equation59} to the boundary conditions 
 \eqref{equation47c}
\begin{align*}
   e_x^w(0,t)=& e^{\tilde{w}}_{x}(0,t)-\int_0^1 k^b_x(0,y)  e^{\tilde{w}}(y,t) dy \nonumber \\
   & +   k^b(0,0) e^{\tilde{w}}(0,t) \\
   =&0,
\end{align*}
where the previous step was obtained by  using  \eqref{equation58b} and $e^{\tilde{w}}(0,t) = e^w(0,t)$. Thus we obtain the  boundary condition in \eqref{equation57c} at $x=0$. Similarly, 
\begin{align*}
   e_x^w(1,t)=& e^{\tilde{w}}_{x}(1,t)-\int_1^1 k^b_x(1,y)  e^{\tilde{w}}(y,t) dy \nonumber\\
   &+ k^b(1,1) e^{\tilde{w}}(1,t) \\
   =& k^b(1,1) e^{\tilde{w}}(1,t)= k^b(1,1) e^w(1,t).  
\end{align*}
If $\eta_3= -k^b(1,1)$ then we obtain the  boundary condition in \eqref{equation57c} at $x=1$.  We perform similar calculations  on the elliptic equation  \eqref{equation47c}. First, we take the spatial derivative of \eqref{equation60},
\begin{align}
%&e_x^v(x,t)=e^{\tilde{v}}_{x}(x,t)-\int_x^1 k_{2x}(x,y) e^{\tilde{v}}(y,t) dy \nonumber \\
%&+ k_2(x,x) e^{\tilde{v}}(x,t) \label{first derivative for ev wv1}\\
          &e^v_{xx}(x,t)= e^{\tilde{v}}_{xx}(x,t)-\int_x^1 k_{2xx}(x,y) e^{\tilde{v}}(y,t)dy   \nonumber\\
          &+ k_{2x}(x,x) e^{\tilde{v}}(x,t) + \frac{d}{dx} k_2(x,x) e^{\tilde{v}}(x,t) \nonumber \\
          &+ k_2(x,x) e^{\tilde{v}}_x(x,t).  \label{equation66}
\end{align}
Subbing \eqref{equation66} in the right-hand-side of elliptic equation  \eqref{equation57b},
\begin{align}
& \text{(R.H.S)}_2=e^{\tilde{v}}_{xx}(x,t)-\int_x^1 k_{2xx}(x,y) e^{\tilde{v}}(y,t)dy \nonumber \\
& + k_{2x}(x,x) e^{\tilde{v}}(x,t)  + \frac{d}{dx} k_2(x,x) e^{\tilde{v}}(x,t)+ k_2(x,x) \nonumber \\ 
& \times e^{\tilde{v}}_x(x,t)   - \gamma e^{\tilde{v}}(x,t)+ \gamma \int_x^1 k_2(x,y) e^{\tilde{v}}(y,t) dy \nonumber \\
& + \beta e^{\tilde{w}}(x,t) - \beta \int_x^1 k_1(x,y) e^{\tilde{w}}(y,t) dy   - \eta_2(x) e^{\tilde{v}}(1,t)\label{equation67}\\
     & \text{(L.H.S)}_2=0. \label{left handside}
\end{align}
Rewriting the last term of \eqref{equation67} as follows 
\begin{align*}
   & \beta \int_x^1  k_1(x,y) e^{\tilde{w}}(y,t) dy  =   - \int_x^1  k_1(x,y) e^{\tilde{v}}_{yy}(y,t) dy \nonumber \\
   &+ (o_2 +\gamma) \int_x^1  k_1(x,y) e^{\tilde{v}}(y,t) dy,
\end{align*}
which can be obtained by referring to the elliptic equation of \eqref{equation47b}, then  \eqref{equation67}  gives
\begin{align}
& \text{(R.H.S)}_2=e^{\tilde{v}}_{xx}(x,t)-\int_x^1 k_{2xx}(x,y) e^{\tilde{v}}(y,t)dy  \nonumber \\
& + k_{2x}(x,x) e^{\tilde{v}}(x,t)  + \frac{d}{dx} k_2(x,x) e^{\tilde{v}}(x,t)+ k_2(x,x)   \nonumber \\
     & \times e^{\tilde{v}}_x(x,t)- \gamma e^{\tilde{v}}(x,t)+ \gamma \int_x^1 k_2(x,y) e^{\tilde{v}}(y,t) dy \nonumber \\ 
     & + \beta e^{\tilde{w}}(x,t) + \int_x^1  k_1(x,y) e^{\tilde{v}}_{yy}(y,t) dy   - \eta_2(x) e^{\tilde{v}}(1,t) \nonumber \\
     & -(o_2 +\gamma) \int_x^1  k_1(x,y) e^{\tilde{v}}(y,t) dy.\label{equation68}
\end{align}
 Since $k_1(x,y)= k_2(x,y)=k^b(x,y)$ , then \eqref{equation68} leads to 
\begin{align*}
    &\text{(R.H.S)}_2=e^{\tilde{v}}_{xx}(x,t)-\int_x^1 k^b_{xx}(x,y) e^{\tilde{v}}(y,t)dy  \nonumber \\ 
    &  + k^b_x(x,x)   e^{\tilde{v}}(x,t)  + \frac{d}{dx} k^b(x,x) e^{\tilde{v}}(x,t)  + k^b(x,x)    \nonumber \\
    & \times e^{\tilde{v}}_x(x,t)- \gamma e^{\tilde{v}}(x,t)+ \gamma \int_x^1 k^b(x,y) e^{\tilde{v}}(y,t) dy \nonumber \\
    & + \beta e^{\tilde{w}}(x,t)  -(o_2 +\gamma) \int_x^1  k^b(x,y) e^{\tilde{v}}(y,t) dy  - \eta_2(x)  \nonumber \\
    &  \times e^{\tilde{v}}(1,t) +  e^{\tilde{v}}_x(1,t) k^b (x,1)  -  e^{\tilde{v}}_x(x,t) k^b(x,x) - e^{\tilde{v}}(1,t)  \nonumber \\
    & \times k^b_y(x,1)  +  e^{\tilde{v}}(x,t) k^b_y(x,x)  + \int_x^1  k^b_{yy}(x,y) e^{\tilde{v}}(y,t) dy .
\end{align*}
Adding and subtracting the term $o_2 e^{\tilde{v}}(x,t)$ and  incorporating $e^{\tilde{v}}_x(1,t)=0$,
\begin{align}
     &\text{(R.H.S)}_2= e^{\tilde{v}}_{xx}(x,t) -   (o_2+\gamma) e^{\tilde{v}} (x,t)+  \beta e^{\tilde{w}}(x,t)   \nonumber \\
     &+\int_x^1 [  -k_{xx}(x,y) +  k_{yy}(x,y) - o_2  k(x,y)] e^{\tilde{v}}(y,t)dy \nonumber \\
     &   + (2 \frac{d}{dx} k^b(x,x)  + o_2)e^{\tilde{v}}(x,t)   - ( k^b_y(x,1)  + \eta_2(x)) \nonumber \\
     &  \times e^{\tilde{v}}(1,t).\label{equation53}
\end{align}
Since $k^b(x,y)$ is given by \eqref{equation58} and  $\eta_2(x) = - k^b_y(x,1)$, then referring to \eqref{equation47b} and \eqref{left handside}
\begin{align*}
  \text{(L.H.S)}_2=0=  \text{(R.H.S)}_2.
\end{align*}

Thus the state transformation \eqref{transformationwv1} transforms the elliptic equation \eqref{equation47b} into \eqref{equation57b}. 
We apply the transformation to the boundary conditions 
 \eqref{equation47d}, %. Referring to \eqref{first derivative for ev wv1} , 
\begin{align*}
   e_x^v(0,t)=& e^{\tilde{v}}_{x}(0,t)-\int_0^1 k^b_x(0,y)  e^{\tilde{v}}(y,t) dy \nonumber \\
   &+k^b(0,0) e^{\tilde{v}}(0,t) =0, 
\end{align*}
by means of using  \eqref{equation58b} and that $ e^{\tilde{v}}_x(0,t) =0$ . We obtain the boundary condition at $x=0$ in \eqref{equation57d}.
 Similarly, 
\begin{align*}
   e_x^v(1,t)=& e^{\tilde{v}}_{x}(1,t)-\int_1^1 k^b_x(1,y)  e^{\tilde{v}}(y,t) dy + k^b(1,1) \nonumber \\
   & e^{\tilde{v}}(1,t) \\
   =& k^b(1,1) e^{\tilde{v}}(1,t) = k^b(1,1) e^v(1,t) ,
\end{align*}
where the previous step was obtained by noting that $e^{\tilde{v}}(0,t) = e^v(0,t)$. If $\eta_4= -k^b(1,1)$ then we obtain the second boundary condition in \eqref{equation57d} at $x=1$. The conclusion of the theorem follows. $\hfill \square$

 The next theorem follows from Theorem \ref{Theorem 11wv1}. 
\begin{thm}
   Let $k(x,y)$ be the solution of system \eqref{equation58}. The error dynamics \eqref{equation57} with output injections $\eta_j, \; j = 1, \dots, 4$ defined as given in \eqref{outputinjection2.1}-\eqref{outputinjection2.2}  are exponentially stable if and only if the parameter  $o_2$ satisfies
    \begin{align}
   (o_2+\rho)(o_2+\gamma)> \alpha \beta, \label{critera for c1}
\end{align}
    and  $o_2 + \gamma \neq -(n \pi)^2 $. 
\end{thm}
\noindent \textbf{Proof.} If $o_2$ is given by \eqref{critera for c1} such that $o_2 +\gamma \neq -(n \pi)^2 $, the target system \eqref{equation47} has a unique solution and is exponentially stable due to the criteria for stability of parabolic-elliptic systems established previously in  Corollary 2 (of a previous draft). Finally, the exponential stability of the error dynamics \eqref{equation57} follows  by referring to Theorem \ref{Theorem 11wv1} and using the invertiblity of transformation \eqref{transformationwv1}. This concludes the proof. $\hfill \square$

\subsection{Only $w(1,t)$ is available} \label{subsection2}
 The objective of this subsection is to design an exponentially convergent observer for \eqref{1}-\eqref{4} given only a single measurement $w(1,t)$. We propose the following observer 
\begin{subequations}\label{2newlabelequation69}
\begin{align}
     \hat{w}_t(x,t) =&  \hat{w}_{xx}(x,t) - \rho  \hat{w}(x,t)  + \alpha \hat{v}(x,t)  \nonumber \\
     & +\eta_1(x) [ w(1,t)- \hat{w}(1,t)], \label{2newlabelequation69a}\\
    0 =&  \hat{v}_{xx}(x,t)  - \gamma \hat{v}(x,t) + \beta \hat{w}(x,t) , \label{2newlabelequation69b}\\
   \hat{w}_x(0,t) =& 0, \quad \hat{w}_x(1,t) =u(t)+\eta_2 [ w(1,t)- \hat{w}(1,t)], \label{2newlabelequation69c} \\
   \hat{v}_x (0,t) =&0, \quad  \hat{v}_x(1,t) =0, \label{2newlabelequation69d} 
\end{align} 
\end{subequations}
where  $\eta_1(x)$ and $\eta_2 $ are   output injections to be designed. 
Defining the  states of the error dynamics as in \eqref{error states}, the system describing the observation error  satisfies
\begin{subequations} \label{2newlabelequation74}
    \begin{align}
        e^w_t(x,t) =& e^w_{xx}(x,t) - \rho e^w(x,t) + \alpha e^v(x,t) \nonumber \\
        & - \eta_1(x) e^w(1,t),  \label{2newlabelequation74a} \\
        0=& e^v_{xx} (x,t) - \gamma e^v (x,t) + \beta e^w(x,t), \label{2newlabelequation74b} \\
        e^w_x(0,t)=& 0,  \quad e^w_x(1,t)=-\eta_2 e^w(1,t),  \label{2newlabelequation74c}\\
         e^v_x(0,t)=& 0,  \quad   e^v_x(1,t)=0 \label{2newlabelequation74d}. 
    \end{align}
\end{subequations}
Both of $\eta_1(x)$ and $\eta_2$ have to be chosen so that  exponential stability of  error dynamics is achieved.  Following a backstepping approach, we define the transformation
\begin{align}
e^{\tilde{w}}(x,t)  = e^w(x,t) - \int_0^x k^a(x,y) e^w(y,t) dy , \label{2newlabelequation75-rewrite}
\end{align}
where $k^a(x,y)$ is given by  \eqref{k-kernel} with $c_2$ replaced by $o_2$. The inverse transformation \cite[Chap. 4, section 5]{krstic2008boundary} is
\begin{align}
e^w(x,t)  = e^{\tilde{w}}(x,t) + \int_0^x \ell^a(x,y) e^{\tilde{w}}(y,t) dy, \label{newlabelinverse-w(1)}
\end{align}
where $\ell^a$ satisfies system \eqref{inverse systemab}.

\begin{thm}  \label{2newlabelTheorem17}
If the output injections are 
\begin{align}
&\eta_1(x) = 0, \label{2newlabeloutputinjection eta1}\\ 
 & \eta_2= -k^a(1,1),\label{2newlabeloutputinjection3.1}  
\end{align}
where $k^a(x,y)$ is given in \eqref{k-kernel} with $c_2$ being replaced by $o_2$  then transformation \eqref{2newlabelequation75} converts the error dynamics \eqref{2newlabelequation74} into 
the target system
\begin{subequations}\label{2newlabel22}
\begin{align}
    & e^{\tilde{w}}_t(x,t) = e^{\tilde{w}}_{xx}(x,t) - (o_2+\rho) e^{\tilde{w}}(x,t)  + \alpha e^v(x,t) \nonumber \\
     & - \alpha \int_0^x k^a(x,y)e^v(y,t) dy, \label{2newlabel22a}\\
    &0 =  e^v_{xx}(x,t)  - \gamma e^v(x,t) + \beta e^{\tilde{w}}(x,t)   \nonumber \\
     &  + \beta  \int_0^x \ell^a (x,y) e^{\tilde{w}}(y,t)  dy,  \label{2newlabel22b}\\
   &e^{\tilde{w}}_x(1,t) =-\int_0^1 k^a_x(1,y)e^{\tilde{w}}(y,t)  dy \nonumber \\
   &-\int_0^1 k^a_x(1,y) \int_0^y \ell^a(y,z) e^{\tilde{w}}(z,t) dz dy, \label{2newlabel22c}\\
   & e^{\tilde{w}}_x(0,t) = 0, \quad e^v_x (0,t) =0, \quad e^v_x (1,t) =0.  \label{2newlabel22d}
\end{align}
\end{subequations}
\end{thm}
\noindent \textbf{Proof.}
It will be useful to rewrite \eqref{2newlabelequation75-rewrite} as
\begin{align}
    e^w(x,t) = e^{\tilde{w}}(x,t)  + \int_0^x k^a(x,y) e^w(y,t) dy. \label{2newlabelequation75}
\end{align}
We take the spatial and the time  derivatives of \eqref{2newlabelequation75} we have
\begin{align}
          &e^w_{xx}(x,t) = e^{\tilde{w}}_{xx}(x,t) +\int_0^x k^a_{xx}(x,y) e^w(y,t)dy   \nonumber \\ 
     & + k^a_x(x,x) e^w(x,t) + \frac{d}{dx} k^a(x,x) e^w(x,t) \nonumber \\
     &+ k^a(x,x) e^w_x(x,t),  \label{2newlabel29}
\end{align}
\begin{align}
         &e^w_{t}(x,t)= e^{\tilde{w}}_{t}(x,t) +\int_0^x k^a(x,y) e^w_t(y,t)dy \nonumber \\
      &=   e^{\tilde{w}}_{t}(x,t)  - \rho \int_0^x k^a(x,y) e^w(y,t)dy   + \alpha \int_0^x k^a(x,y) \nonumber \\
      & \times e^v(y,t)dy  + k^a(x,x) e^w_x(x,t) 
      - k^a(x,0) e^w_x(0,t)  \nonumber \\
      & - k^a_y(x,x) e^w(x,t) +k^a_y(x,0) e^w(0,t) +   \int_0^x  k^a_{yy} (x,y) \nonumber \\
      & \times e^w(y,t) dy -  e^w(1,t) \int_0^x k^a(x,y) \eta_1(y) dy .
         \label{2newlabel30}
\end{align}
Substituting \eqref{2newlabel29} and \eqref{2newlabel30} in the parabolic equation \eqref{2newlabelequation74a},
\begin{comment}
\begin{align*}
     &e^{\tilde{w}}_{t}(x,t)  - \rho \int_0^x k^a(x,y) e^w(y,t)dy   + \alpha \int_0^x k^a(x,y) \nonumber \\
     & \times e^v(y,t)dy  + k^a(x,x) e^w_x(x,t) 
       - k^a(x,0) e^w_x(0,t) \nonumber \\
      & - k^a_y(x,x) e^w(x,t) +k^a_y(x,0) e^w(0,t) +   \int_0^x  k^a_{yy} (x,y) \nonumber \\
      & \times e^w(y,t) dy -  e^w(1,t) \int_0^x k^a(x,y) \eta_1(y) dy \\
      &= e^{\tilde{w}}_{xx}(x,t) +\int_0^x k^a_{xx}(x,y) e^w(y,t)dy + k^a_x(x,x) e^w(x,t) \nonumber \\
      & + \frac{d}{dx} k^a(x,x) e^w(x,t)  + k^a(x,x) e^w_x(x,t) - \rho  e^w(x,t) \nonumber \\ 
     & + \alpha e^v(x,t) - \eta_1(x) e^w(1,t). 
\end{align*}
\end{comment}
and using  $e^w_x(0,t) =0$ and  $k^a_y(x,0)=0$, %one can rearrange the terms in the equation above to obtain
\begin{align}
     & e^w_t(x,t)   = e^w_{xx}(x,t) + \alpha e^v(x,t) +\left(k^a_y(x,x) +k^a_x(x,x)\right) \nonumber \\
     & \times e^w(x,t) + \frac{d}{dx} k^a(x,x) e^w(x,t)  - \rho  e^w(x,t)  + k^a(x,x) \nonumber \\
     & \times e^w_x(x,t) - k^a(x,x) e^w_x(x,t) 
      + \int_0^x[k^a_{xx} (x,y) \nonumber \\
      & -k^a_{yy}(x,y) +\rho k^a(x,y) ] e^w(y,t)dy  - \alpha \int_0^x k^a(x,y) \nonumber \\
      & \times e^v(y,t)dy  + e^w(1,t) \int_0^x k^a(x,y) \eta_1(y) dy  \nonumber \\
      &- \eta_1(x) e^w(1,t)  . \label{newlabel01-w(1)}
\end{align}
Adding and subtracting the term $(o_2+\rho) e^w(x,t) $ to the right-hand-side of equation \eqref{newlabel01-w(1)}, and noting that  $k^a(x,y)$ is given by  \eqref{26},
\begin{comment}
%\begin{align*}
%-(o_2+\rho) e^w(y,t)= & - (o_2+\rho)e^{\tilde{w}}(x,t)  \nonumber \\
%&-(o_2+\rho) \int_x^1 k^a(x,y) e^w(y,t) dy,
%\end{align*}
%then
%\begin{align*}
%     &e^{\tilde{w}}_t(x,t)   = e^{\tilde{w}}_{xx}(x,t)  - (o_2+\rho) e^{\tilde{w}}(x,t)  + \alpha e^v(x,t) \nonumber \\
 %    & - \alpha \int_0^x k^a(x,y) e^v(y,t)dy +  \int_0^x [ k^a_{xx} (x,y) - k^a_{yy}(x,y) \nonumber \\
  %   & - o_2 k^a(x,y)] e^w(y,t) dy + (2  \frac{d}{dx} k^a(x,x)  +o_2)   e^w(x,t) \\
%& + e^w(1,t) \int_0^x k^a(x,y) \eta_1(y) dy - \eta_1(x) e^w(1,t) .
%\end{align*}
\end{comment}
%Since $k^a(x,y)$ is given by  \eqref{26}
\begin{align*}
     & e^w_t(x,t)   = e^w_{xx}(x,t)  - \rho e^w(x,t)  + \alpha e^v(x,t) \nonumber \\
     & - \alpha \int_0^x k^a(x,y) e^v(y,t)dy  - e^w(1,t) \int_0^x k^a(x,y)  \nonumber \\
     & \times \eta_1(y) dy  - \eta_1(x) e^w(1,t) .
\end{align*}
If $\eta_1(x)=0$,  we obtain the parabolic equation \eqref{2newlabel22a}.   We now apply transformation \eqref{2newlabelequation75} on the boundary conditions  \eqref{2newlabelequation74c}, using lemma \ref{Wellposedness}
\begin{align*}
      e^{\tilde{w}}_x(0,t) &= e^w_x(0,t) - \int_0^0 k^a_x(1,y) w(y)dy \nonumber \\
      & - k^a(0,0) w(1,t) = 0,
\end{align*}
by virtue of referring to the boundary conditions of system \eqref{26}.
\begin{align*}
    &e^{\tilde{w}}_x(1,t) = e^w_x(1,t) - \int_0^1 k^a_x(1,y) e^w(y,t)dy \\
    &- k^a(1,1) e^w(1,t) \\
    =&-(\eta_2+k^a(1,1))e^w(1,t) - \int_0^1 k^a_x(1,y) e^w(y,t)dy\\
    =& - \int_0^1 k^a_x(1,y) e^w(y,t)dy .\\
   =& -\int_0^1 k^a_x(1,y)e^{\tilde{w}}(y,t)  dy-\int_0^1 k^a_x(1,y) \int_0^y \ell^a(y,z) \nonumber \\
   & \times e^{\tilde{w}}(z,t) dz dy.
\end{align*} The previous  equation holds true via using \eqref{2newlabeloutputinjection3.1} and using the inverse transformation \eqref{newlabelinverse-w(1)}. The elliptic equation  \eqref{2newlabel22b} can be obtained via using the inverse transformation \eqref{newlabelinverse-w(1)}.
$\hfill \square$

\begin{thm}\label{2newlabelTheorem18}
    The target system \eqref{2newlabel22} is exponentially stable if $o_2$ is chosen such that
    \begin{align}
        o_2+\rho&> \frac{|\alpha||\beta|}{\gamma}   (1+\|k^a\|) (1+\|\ell^a\|)\nonumber \\
        &+ \frac{(1+\|\ell^a\|)^2\|k^a_x(1,y)\|^2+2}{2}. \label{1c_3-w(1)}
    \end{align}
\end{thm}
\noindent \textbf{Proof.}  We define the Lyapunov function candidate,
\begin{align*}
    V(t) =& \frac{1}{2} \int_0^1 (e^{\tilde{w}}(x,t))^2dx = \frac{1}{2} \| e^{\tilde{w}}(x,t) \|^2  .
\end{align*}
Taking the time derivative of $V(t)$,
\begin{align}
    & \dot{V}(t) =\int_0^1 e^{\tilde{w}}(x,t) e^{\tilde{w}}_t(x,t) dx  \nonumber \\
    = & \int_0^1 e^{\tilde{w}}(x,t) e^{\tilde{w}}_{xx}(x,t)dx -(o_2+\rho) \int_0^1 (e^{\tilde{w}}(x,t))^2 dx \nonumber \\
    & + \alpha \int_0^1 e^{\tilde{w}}(x,t) e^v(x,t) dx - \alpha \int_0^1 e^{\tilde{w}}(x,t)\int_0^x k^a(x,y) \nonumber \\
    & \times e^v(y,t) dydx . \label{2newlabel31}
\end{align}
Integrating the term $\int_0^1 e^{\tilde{w}}(x,t) e^{\tilde{w}}_{xx}(x,t)dx$ by parts, and using the boundary conditions  \eqref{2newlabel22c}-\eqref{2newlabel22d}
\begin{align}
    &\int_0^1  e^{\tilde{w}}(x,t) e^{\tilde{w}}_{xx}(x,t)dx=e^{\tilde{w}}(1,t) e^{\tilde{w}}_{x}(1,t) \nonumber \\
    & -e^{\tilde{w}}(0,t) e^{\tilde{w}}_{x}(0,t)- \| e^{\tilde{w}}_x\|^2 \nonumber \\
    &= e^{\tilde{w}}(1,t) e^{\tilde{w}}_{x}(1,t)- \| e^{\tilde{w}}_x\|^2 \nonumber \\
    &= -\int_0^1 k^a_x(1,y) e^w(y,t)dy \,e^{\tilde{w}}(1,t) - \| e^{\tilde{w}}_x\|^2. \label{2newlabeleq17}
\end{align}
To bound the term $-\int_0^1 k^a_x(1,y) e^w(y,t)dy \,e^{\tilde{w}}(1,t)$ in \eqref{2newlabeleq17}, we use Cauchy-Schwartz, 
\begin{align*}
   &- \int_0^1 k^a_x(1,y) e^w(y,t)dy e^{\tilde{w}}(1,t)  \nonumber \\
   &\leq \|\int_0^1 k^a_x(1,y) e^w(y,t)dy \|  \max_{x\in[0,1]} e^{\tilde{w}}(x,t) \nonumber \\ 
    & \leq \| k^a_x(1,y)\| \|e^w\| \| e^{\tilde{w}} \|_\infty .
\end{align*}
Invoking Agmon's inequality \cite{triggiani1975stabilizability} on the right-hand-side of the previous inequality leads to
\begin{align*}
    & - \int_0^1 k^a_x(1,y) e^w(y,t)dy e^{\tilde{w}}(1,t) \nonumber \\
    &\leq \| k^a_x(1,y)\| \|e^w\| \| e^{\tilde{w}}\|^{1/2} \| e^{\tilde{w}}\|_{H^1}^{1/2}.
\end{align*}
Using Young's inequality on the right-hand-side of the previous inequality,
\begin{align*}
  & - \int_0^1 k^a_x(1,y) e^w(y,t)dy e^{\tilde{w}}(1,t) \nonumber \\
  & \leq \frac{\| k^a_x(1,y)\|^2}{2} \|e^w\|^2 + \frac{1}{2} (\| e^{\tilde{w}}\| \| e^{\tilde{w}}\|_{H^1}) \nonumber \\
    & \leq \frac{\| k^a_x(1,y)\|^2}{2} \|e^w\|^2 + \frac{1}{2} \| e^{\tilde{w}}\|^2 +  \frac{1}{2} \| e^{\tilde{w}}\| \| e^{\tilde{w}}_x\| \nonumber .
   \end{align*}
Estimating the term $\frac{1}{2} \| e^{\tilde{w}}\| \| e^{\tilde{w}}_x\|$ on the right-hand-side of the previous inequality using Young's inequality, then
\begin{align}
    &- \int_0^1 k^a_x(1,y) e^w(y,t)dy e^{\tilde{w}}(1,t)     \leq \frac{\| k^a_x(1,y)\|^2}{2} \|e^w\|^2 \nonumber \\
    &+  \|e^{\tilde{w}}\|^2 + \frac{1}{4} \| e^{\tilde{w}}_x\|^{2} .
\end{align}
Referring to the inverse transformation \eqref{newlabelinverse-w(1)},
\begin{align}
& - \int_0^1 k^a_x(1,y) e^w(y,t)dy e^{\tilde{w}}(1,t)     \nonumber \\
&\leq \frac{(1+\|\ell^a\|)^2\|k^a_x(1,y)\|^2+2}{2}  \|e^{\tilde{w}}\|^2 + \frac{1}{4} \| e^{\tilde{w}}_x\|^{2} . \label{2newlabel15}
\end{align}
Combining \eqref{2newlabel15} and \eqref{2newlabeleq17}, we get
\begin{align}
    &\int_0^1  e^{\tilde{w}}(x,t) e^{\tilde{w}}_{xx}(x,t)dx  \nonumber \\
    & \leq \frac{(1+\|\ell^a\|)^2\|k^a_x(1,y)\|^2+2}{2}  \|e^{\tilde{w}}\|^2 -(1- \frac{1}{2}) \| e^{\tilde{w}}_x\|^{2} \nonumber \\
   & \leq \frac{(1+\|\ell^a\|)^2\|k^a_x(1,y)\|^2+2}{2}  \|e^{\tilde{w}}\|^2. \label{2newlabelequation19}
\end{align}
Bounding the  terms on the right-hand side of \eqref{2newlabel31} using \eqref{2newlabelequation19}, Cauchy-Schwartz inequality and lemma \ref{lemma}, we arrive to 
\begin{align}
   \dot{V}(t) \leq&   -(o_2+\rho- \frac{|\alpha||\beta|}{\gamma}   (1+\|k^a\|) (1+\|\ell^a\|) \nonumber \\
   & - \frac{(1+\|\ell^a\|)^2\|k^a_x(1,y)\|^2+2}{2} ) \| e^{\tilde{w}}\|^2.\label{2newlabel34}
\end{align}

Setting
\begin{align*}
o_3=& o_2+\rho- \frac{|\alpha||\beta|}{\gamma}   (1+\|k^a\|) (1+\|\ell^a\|) \nonumber \\
& - \frac{(1+\|\ell^a\|)^2\|k^a_x(1,y)\|^2+2}{2},
\end{align*}
then inequality \eqref{2newlabel34} implies that $ V(t) \leq  e^{-2 o_3  t} V(0)$. If the parameter $o_2$ is chosen such that \eqref{1c_3-w(1)} is satisfied, then $V(t)$ decays exponentially as $t \rightarrow \infty$. Thus, $\| e^{\tilde{w}}(x,t)\|$ decays exponentially. Referring to the elliptic equation of system \eqref{2newlabelequation74} and recalling lemma \ref{lemma},  the state $e^v(x,t)$ is asymptotically stable. Exponential stability of system \eqref{2newlabel22}  follows by noting that $\partial_{xx} -\gamma I$ is boundedly invertible, and using  a similar argument as the one given in the last part in the proof of Theorem \ref{thm-stab}. $\hfill \square$

The following lemma, which establishes a bound on $\|k^a_x(1,y)\|$, will be needed.

\begin{lem} \label{bounds on k_x1}
    Consider system \eqref{26} with $c_2$ being replaced by $o_2$. The $L^2$-norm of  $k^a_x(1,y)$ is bounded by
\begin{align}
     & \| k^a_x(1,y) \| \leq  \frac{o_2}{2}(1+\frac{o_2}{2} ) e^{\frac{o_2 }{4}} \left( \sqrt\frac{\pi}{2o_2}
 erf(\sqrt{\frac{o_2}{2}}) \right)^{\frac{1}{2}} . \label{kx1bound}
 \end{align}
\end{lem}
 \noindent \textbf{Proof.} %%%%%%%%%%%%%%%%%%%%%%%%%%%%%
 The relation \eqref{kx1bound} can be shown by noting that the solution of system \eqref{26} is 
 \begin{align*}
     k^a(x,y)= -o_2 x \frac{I_1(\sqrt{o_2 (x^2-y^2)})}{\sqrt{o_2 (x^2-y^2)}}.
 \end{align*}
After straightforward  mathematical steps, we arrive to
\begin{align}
     k^a_x(x,y)=& -o_2  \frac{I_1(\sqrt{o_2 (x^2-y^2)})}{\sqrt{o_2(x^2-y^2)}} \nonumber \\
     & - o_2 x \frac{I_2(\sqrt{o_2 (x^2-y^2)}}{(x^2-y^2)} \label{k_x-w(1)},
 \end{align}
 where we have used that $\frac{d}{dx} I_1(x)= \frac{I_1(x)}{x}+ I_2(x)$. Setting $z=\sqrt{o_2 (x^2-y^2)}$ and using the definition of Bessel function, \eqref{k_x-w(1)} can be written as

 \begin{align}
     &k^a_x(x,y)= -o_2  \frac{I_1(z)}{z} -o_2^2 x \frac{I_2(z)}{z^2} \nonumber \\
     &=- \frac{o_2}{z} \sum_{m=0}^{\infty} \left( \frac{z}{2}\right)^{2m+1} \frac{1}{m! m+1!} 
     - \frac{o_2^2}{z^2} x  \sum_{m=0}^{\infty} \left( \frac{z}{2}\right)^{2m+2} \nonumber \\
     &\times \frac{1}{m! m+2!} .  \label{new-equaion3}
     %&=   -\frac{o_2}{2}  \sum_{m=0}^{\infty} \left( \frac{z^2}{4}\right)^{m} \frac{1}{m! m+1!} 
     %-  \frac{o_2^2}{4} x  \sum_{m=0}^{\infty} \left( \frac{z^2}{4}\right)^{m}  \nonumber \\
    % &\times\frac{1}{m! m+2!} \nonumber \\
    %&\leq    \frac{o_2}{2}  \sum_{m=0}^{\infty} \left( \frac{z^2}{4}\right)^{m} \frac{1}{m! m+1!} 
    % +\frac{o_2^2}{4} x  \sum_{m=0}^{\infty} \left( \frac{z^2}{4}\right)^{m}   \nonumber \\
    % &\times \frac{1}{m! m+2!} \nonumber \\
 \end{align}
 To find a bound on the induced $L_2$-  norm of $k^a_x(1,y)$, with $x=1$ the variable $z$ becomes  $z=\sqrt{o_2 (1-y^2)}$ where $0<y<1$. Equation \eqref{new-equaion3} leads to
 \begin{align} 
     & \| k^a_x(1,y) \| \leq  \frac{o_2}{2}  \sum_{m=0}^{\infty} \frac{(z^2/4)^m}{m!}
     +  \frac{o_2^2}{4} x \sum_{m=0}^{\infty} \frac{(z^2/4)^m}{m!}  \nonumber \\
     &\leq \frac{o_2}{2} \|e^{\frac{z^2}{4}} \| + \frac{o_2^2}{4} \|e^{\frac{z^2}{4}} \| =\frac{o_2}{2}(1+\frac{o_2}{2} ) \|e^{\frac{z^2}{4}} \| \nonumber \\
     & \leq \frac{o_2}{2}(1+\frac{o_2}{2} )e^{\frac{o_2 }{4}} \|e^{\frac{-o_2 y^2}{4}} \|.
       \label{k_x 2-w(1)}
 \end{align}
  Since  $ \;erf (x) = \frac{2}{\sqrt{\pi}} \int_0^x e^{-\xi^2} d \xi$, inequality \eqref{k_x 2-w(1)} leads to \eqref{kx1bound}.
 \hfill $\square$
 
Using lemma \ref{bounds on k, l} and lemma \ref{bounds on k_x1}, the next corollary to Theorem \ref{2newlabelTheorem18} is now immediate.

\begin{cor}
The observation error dynamics \eqref{2newlabel22}  is exponentially stable if
\begin{align}
   & o_2+\rho  >  \left(\frac{|\alpha||\beta|}{\gamma} + \frac{\|k^a_x(1,y)\|^2}{2}\right)  \nonumber \\
   & \times \left[ 1+\sqrt{\frac{o_2\pi}{2}} \; \left( erfi (\sqrt{\frac{o_2}{2}} ) \right)^{\frac{1}{2}} \left( erf (\sqrt{\frac{o_2}{2}} ) \right)^{\frac{1}{2}}  \right]^2 +1 .
   \label{c3w(1)} 
\end{align} 
\end{cor}

The following result now follows from Theorem \ref{2newlabelTheorem17} and Theorem \ref{2newlabelTheorem18}. 
\begin{thm}
   Let $k^a(x,y)$ be the solution of system \eqref{26}. The error dynamics \eqref{2newlabelequation74} with output injection $\eta_2= \frac{1}{2} o_2 $ defined as given in \eqref{2newlabeloutputinjection3.1}  are exponentially stable if  $o_2$ is given by \eqref{1c_3-w(1)}.
\end{thm}

\begin{figure}[H]
\begin{center}
\textbf{Illustration of the restriction \eqref{c3w(1)}  on $o_2$  }\par\medskip
{\includegraphics[scale=0.3]{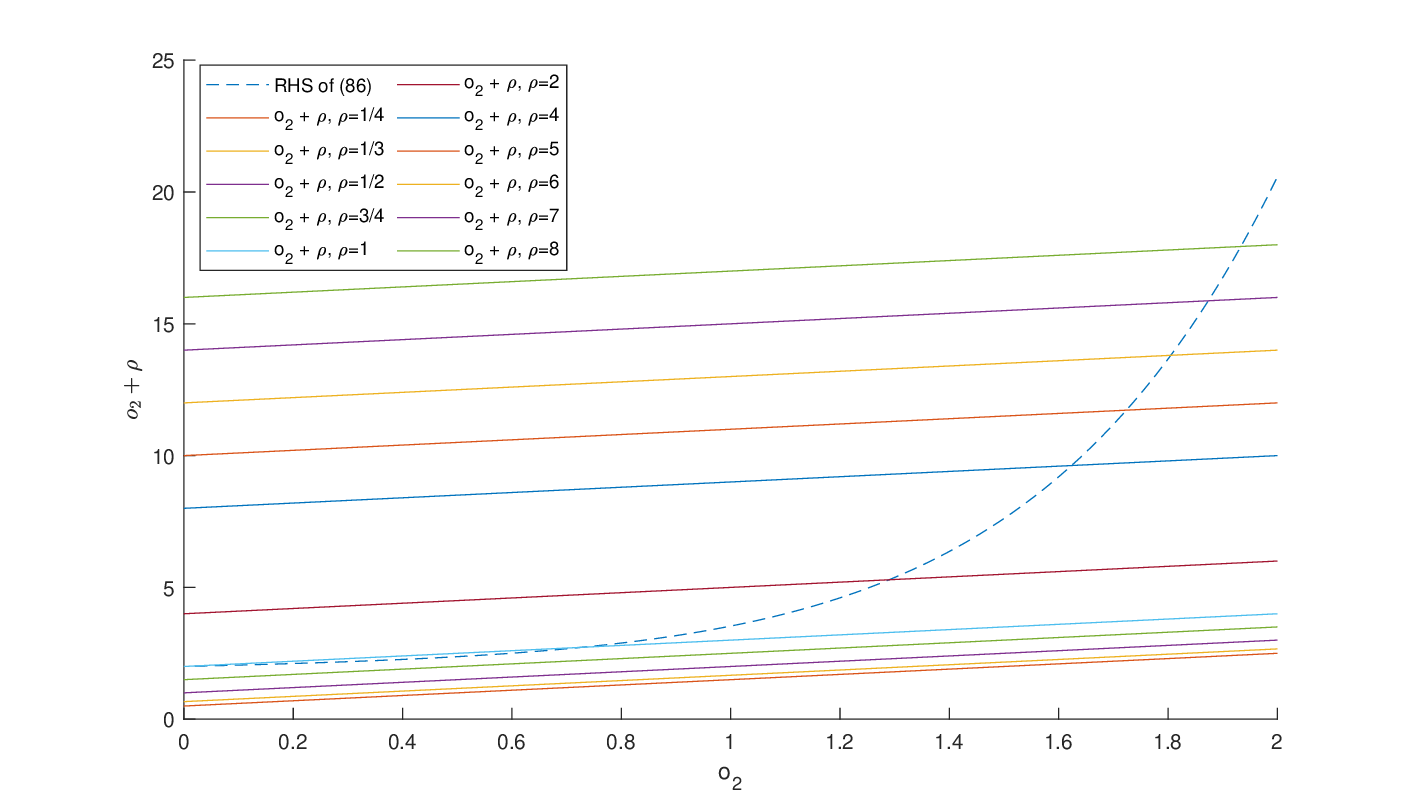} }
\end{center}
\caption{A comparison between the right-hand-side of inequality  \eqref{c3w(1)} as a function of $o_2$ against several straight lines  $o_2+\rho $ for different values of $\rho$, while the other parameters are fixed as $\beta=\gamma=1$, $\alpha=0.5$. The right-hand-side of \eqref{c3w(1)} is described using a dashed line(- - -). For any $\rho$, the error dynamics \eqref{2newlabel22} is exponentially stable if $o_2$ is such that  the dashed line(- - -) is beneath the straight line $o_2 + \rho$. This illustrates the constraints  associated with  condition \eqref{c3w(1)}.}
\label{c2+rho-comparison-single measurement}
\end{figure}

Condition \eqref{c3w(1)} for stability of the observation error dynamics, imposes restrictions on the permissible choices of system parameters.  This observation is demonstrated in Figure\ref{c2+rho-comparison-single measurement}, where we present a  comparison between the right-hand-side of inequality \eqref{c3w(1)} as a function of $o_2$, and several straight lines $o_2 + \rho$ varying with different $\rho$. The other parameters are  fixed as $\beta = \gamma=1, \; \alpha=0.5$. The dashed line in \ref{c2+rho-comparison-single measurement} represents the right-hand-side of \eqref{c3w(1)}. The observation error dynamics \eqref{2newlabelequation74} is exponentially stable if the values of $o_2$ are such that the dashed line in Figure\ref{c2+rho-comparison-single measurement} is beneath the straight lines, for different $\rho$.

\subsection{Numerical simulations}
We conducted numerical simulations for the dynamics of  both the coupled system \eqref{1}-\eqref{4} and the state observer  \eqref{equation56} in the situation where  two measurements, $w(1,t)$ and $v(1,t),$ are available. The simulations were performed  using COMSOL Multiphysic software %As for the  numerical simulations presented in subsection \eqref{numerical simulatino1},   
using linear splines  to approximate the coupled equations by a system of DAEs. The spatial interval was  divided into $27$ subintervals. Time was discretized by  a time-stepping algorithm called  generalized alpha with time-step= 0.1.

\begin{figure}[H]
\begin{center}
\textbf{True and estimated states at $x = 0.56$ using observer \eqref{equation56} }\par\medskip
\subfloat[Comparison between $w(0.5,t)$ and $\hat{w}(0.5,t)$]{\includegraphics[scale=0.25]{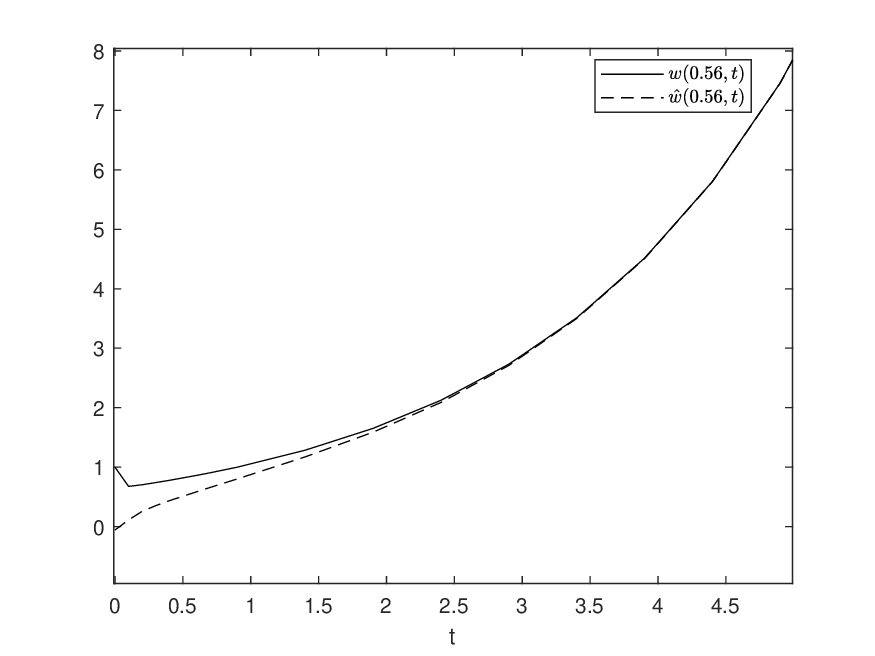}}
\hfill
\subfloat[Comparison between $v(0.5,t)$ and $\hat{v}(0.5,t)$]{\includegraphics[scale=0.25]{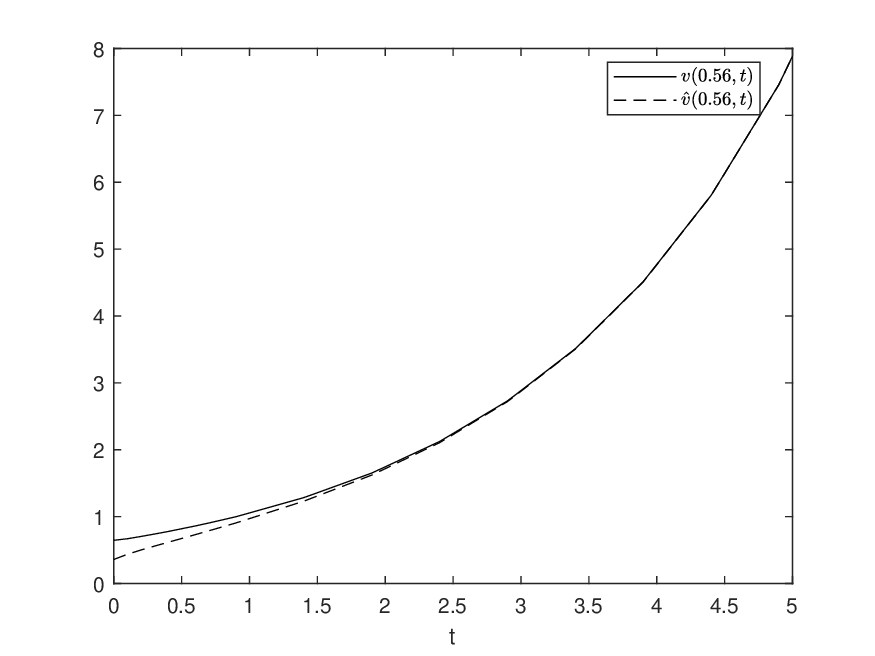}}
\end{center}
\caption{  A comparison between the states of the coupled system \eqref{1}-\eqref{4} versus the estimated states using observer \eqref{equation56} at $x=0.56$. System parameters are $\gamma=1$, $\rho=0.5$,  $\alpha = 1$,  $\beta= 1$ , $o_2=5$ with initial conditions  $w_0=\sin(\pi x)$ and $\hat{w}_0=\cos(\pi x)$.  }
\label{obs-two-meas-midpoint}
\end{figure}

\begin{figure}[H]
\begin{center}
\textbf{$L_2$-norm of the error dynamics \eqref{equation57} using two measurements  $w(1,t), \; v(1,t)$ }\par\medskip
\subfloat[ $\| w(\cdot,t)- \hat{w}(\cdot,t)\|$  ]{\includegraphics[scale=0.25]{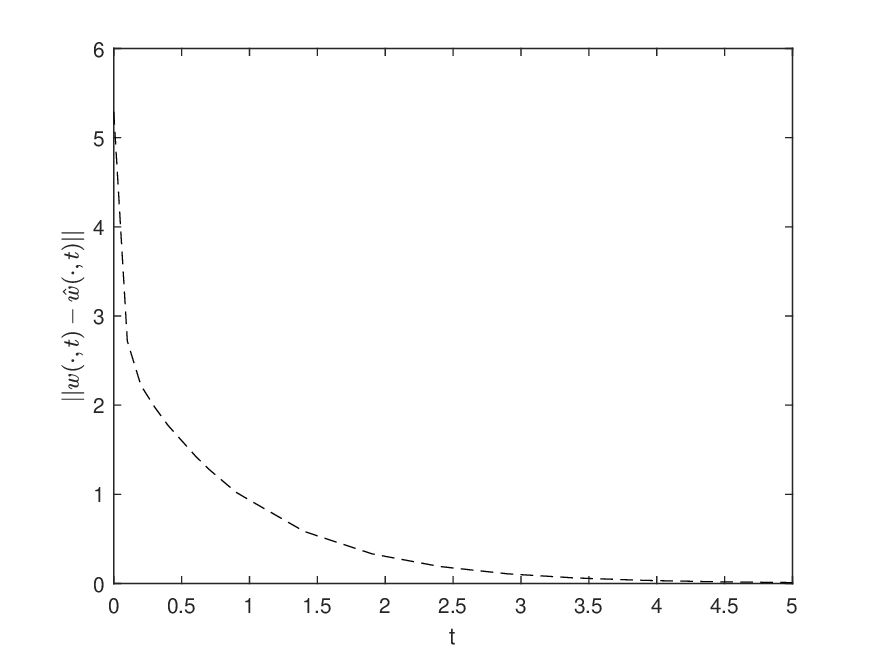}}
\hfill
\subfloat[ $\| v(\cdot,t)- \hat{v}(\cdot,t)\|$]{\includegraphics[scale=0.25]{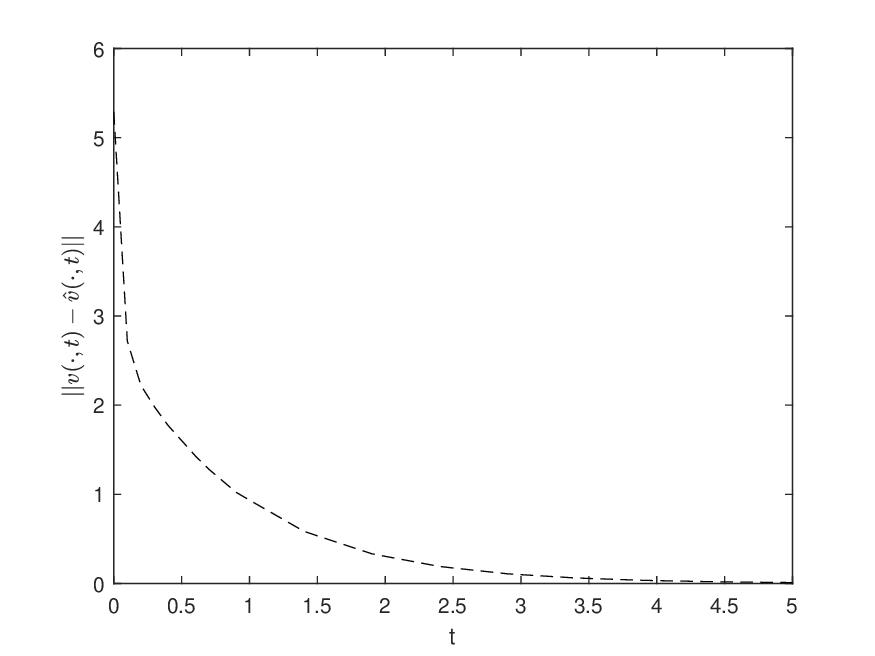}}
\end{center}
\caption{  $L_2$-norms of the error dynamics of the estimates of the parabolic and
elliptic states, using two measurements $w(1,t), \; v(1,t).$  System parameters are $\gamma=1$, $\rho=0.5$,  $\alpha = 1$,  $\beta= 1$ , $o_2=5$ with initial conditions  $w_0=\sin(\pi x)$ and $\hat{w}_0=\cos(\pi x)$.  The estimation error tends to $0$ as $t \to \infty.$}
\label{obs-two-meas-error-dynamics}
\end{figure}

Observer designs were done for system \eqref{1}-\eqref{4} with $u(t) \equiv 0$. The chosen parameters were  $\gamma=1$, $\rho=0.5$,  $\alpha = 1 $ and $\beta= 1 .$  With these parameters, the system is unstable.
%; this is  illustrated in Figure\ref{original system with parameters associated with two measurments-w(1)v(1)}.  
With $o_2= 5,$  the sufficient condition   \eqref{critera for c1} for the error dynamics to be exponentially stable is satisfied. % Figure\ref{observer two measurments-w(1)v(1)} illustrates the observer  with initial condition $\hat{w}_0=\cos(\pi x)$. 
In Figure\ref{obs-two-meas-midpoint} the true and estimated states at $x=0.56$ are shown.
Figure\ref{obs-two-meas-error-dynamics} illustrates the $L_2-$norm error dynamics, which converge to zero as predicted by theory.

%%%%%%%%%%%%% SIMULATIONS  Observer dynamics when w(1) and v(1) are available, alpha=beta=1.5, gamma=1, rho=1.75, c1=5

%% original system with parameters associated with two measurments-w(1)v(1) figure moved to after end document, also L2 norms of both w and w hat etc
\begin{comment}
\begin{figure}[H]
\begin{center}
\textbf{Dynamics of observer \eqref{equation56} 
 with two measurements, $w(1,t)$ and $v(1,t)$}\par\medskip
\subfloat[$\hat{w}(x,t)$]{\includegraphics[scale=0.25]{hat w(x)- when w(1) v(1) are available.eps}}
\hfill
\subfloat[$\hat{v}(x,t)$]{\includegraphics[scale=0.25]{hat v- when w(1) v(1) are available.eps}}
\end{center}
\caption{A 3D landscape of the dynamics of observer \eqref{equation56} with measurements $w(1,t)$ and $v(1,t)$  and initial condition  $\hat{w}_0=\cos(\pi x).$ The system parameters are   $\gamma=1$, $\rho=0.5$,  $\alpha = 1$,  $\beta= 1 $ so the system is unstable. 
%Also, $o_2=5$ so the system meets the condition established in \eqref{critera for c1} for the stability of the  observation error. Observer \eqref{equation56} mimics the unstable dynamics of the coupled system \eqref{1}-\eqref{4} given in Figure\ref{original system with parameters associated with two measurments-w(1)v(1)}.
}
\label{observer two measurments-w(1)v(1)}
\end{figure}
\end{comment}

%

%%%%%%%%%%%%%%%%%%%%%%%%%%%%%%%%%%%%%%%%%%%%%%%%%%%%%%%%%%%%%%%%%%%%% '

\begin{figure}[H]
\begin{center}
\textbf{True and estimated states at $x = 0.56$ using observer \eqref{2newlabelequation69} }\par\medskip
\subfloat[Comparison between $w(0.5,t)$ and $\hat{w}(0.5,t)$]{\includegraphics[scale=0.25]{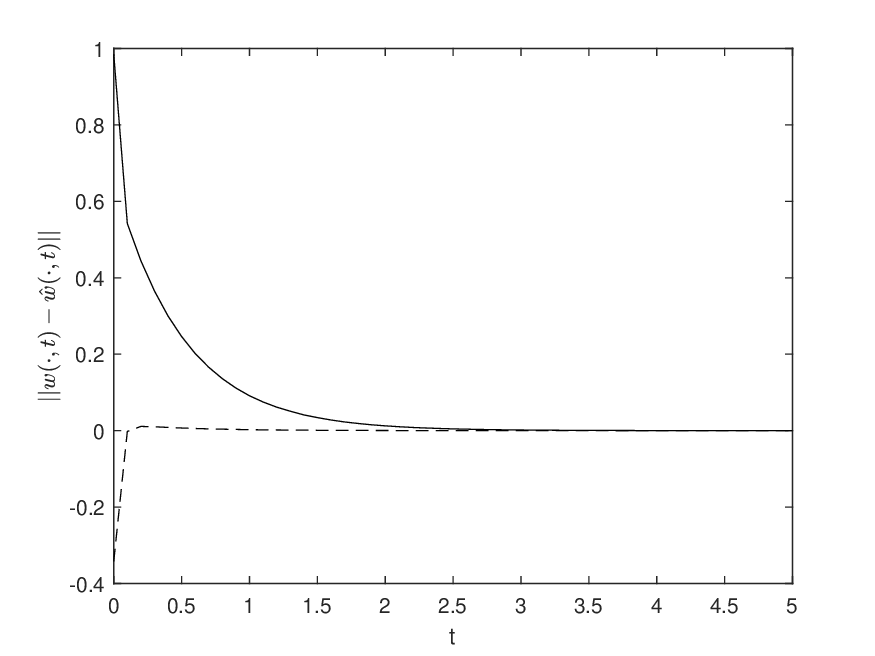}}
\hfill
\subfloat[Comparison between $v(0.5,t)$ and $\hat{v}(0.5,t)$]{\includegraphics[scale=0.25]{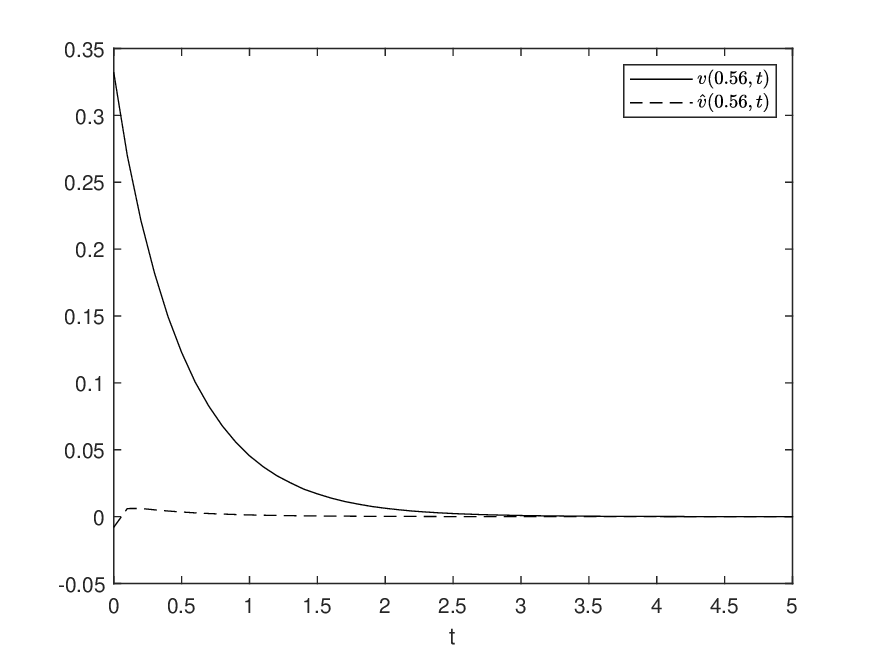}}
\end{center}
\caption{  A comparison between the states of the coupled system \eqref{1}-\eqref{4} versus the estimated states using observer \eqref{2newlabelequation69} at $x=0.56$. Here $\gamma=1$, $\rho=1$,  $\alpha = 0.5$,  $\beta= 0.5$ , $o_2=c_2=0.5$ with $w_0=\sin(\pi x)$ and $\hat{w}_0=\sin(2\pi x)$.  }
\label{observer-states-midpoint-single measurement}
\end{figure}

\begin{figure}[H]
\begin{center}
\textbf{$L_2$-norm of the error dynamics \eqref{2newlabelequation74}, using one measurements  $w(1,t)$ }\par\medskip
\subfloat[   $\|w(x,t)- \hat{w}(x,t)\|$ ]{\includegraphics[scale=0.25]{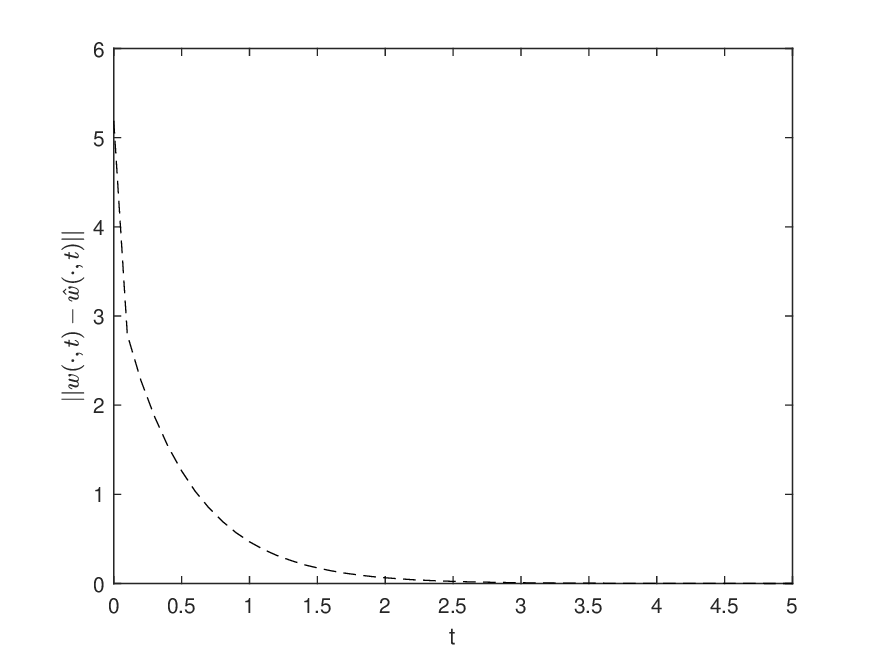}}
\hfill
\subfloat[  $\|v(x,t)- \hat{v}(x,t)\|$]{\includegraphics[scale=0.25]{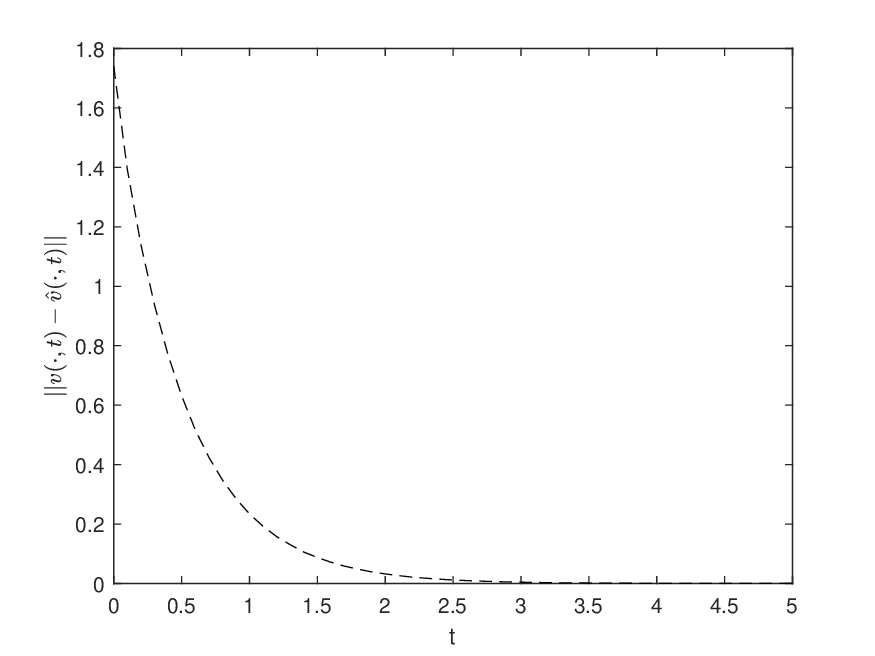}}
\end{center}
\caption{  $L_2$-norms of the error dynamics of the estimates of the parabolic and
elliptic states, using a single measurement $w(1,t).$  Here  $\gamma=1$, $\rho=1$,  $\alpha = 0.5$,  $\beta= 0.5 $ and $o_2=c_2=0.5$  with   $w_0=\sin(\pi x)$ and $\hat{w}_0=\sin(2\pi x)$.  The estimation error tends to $0$ as $t \to \infty.$}
\label{Comparison-observer design with w(1)}
\end{figure}

Numerical simulations were also conducted to study  the observer  \eqref{2newlabelequation69} when a single measurement $w(1,t)$ available. %The results are shown in Figure\ref{original system one measurments-w(1)} and Figure\ref{observer one measurment-w(1)} for the plant and the  observer, respectively. 
The simulations were carried out using the parameter values $\gamma=1$, $\rho=1$, $\alpha = 0.5$, $\beta= 0.5$. With these parameter values, the system is stable.  Also, we set $o_2=0.5$ so the stability condition for the observation error dynamics (i.e. \eqref{c3w(1)}) is satisfied.  The control input $u(t)$ was set as stated in \eqref{control signal stab}, with control gain $c_2=0.5$. The  initial conditions were $w_0=\sin(\pi x)$ and $\hat{w}_0=\sin(2 \pi x).$ 
 The true and estimated states at $x = 0.56$ are given in Figure\ref{observer-states-midpoint-single measurement}.  The  $L_2$-norms of the error dynamics are presented in Figure\ref{Comparison-observer design with w(1)}.

%%%
%%

%%%%%%%%%%%%%%%%%%%%%%%%%%%%%%%
\section{Output feedback} \label{section4}
In general, the full state is not available for control. Output feedback is based on using only the available measurements  to stabilize the system. A common approach to output feedback is to combine a stabilizing state feedback with an observer. The estimated state from the observer is used to replace the state feedback $Kz$ by $K \hat z$ where $z $ is the true state and $\hat z$ the estimated state. 

In the situation considered here, if there are two measurements, this leads to the output feedback controller consisting of the observer \eqref{equation56} combined with the state feedback 
\begin{equation}
u(t)= \int_0^1 k^a_x(1,y) \hat{w}(y,t)dy+k^a(1,1) \hat{w}(1,t) \label{cs}
\end{equation}
where $k^a(x,y)$ is the solution of system \eqref{26} with  $c_2$ satisfying the bound \eqref{bound on c1}. 
Since  the original system  \eqref{1}-\eqref{4} is a well-posed control system (Theorem \ref{1})and  also the observer combined with the state feedback is a well-posed system, the following result follows immediately from the results in \cite[section 3]{morris1994state}.
\begin{thm}
  The closed-loop system consisting of \eqref{1}-\eqref{4}, together with the observer dynamics \eqref{equation56} and control input  \eqref{cs}, is well-posed and exponentially stable.    
\end{thm}

Note that when using output feedback,  the parameters  have to satisfy both of the stability condition associated with the control problem \eqref{bound on c1}, and also the bound of the stability of the observation error  \eqref{critera for c1}.

\subsection{Numerical simulations}
Numerical simulations, again using a finite-element approximation in the COMSOL Multiphysics software, were performed to study the solutions of system (\ref{1})-(\ref{4}) with output feedback \eqref{cs}. The parameter values were $\gamma=\frac{1}{4}$, $\rho=\frac{1}{3}$, $\alpha=\frac{1}{4}$, and $\beta=\frac{1}{2}$. The system's initial condition was  $w_0=\sin(\pi x)$.With these parameters, the uncontrolled system is unstable. %; see subsection \eqref{numerical simulatino1} and in particular Figure\ref{fig1}.

\begin{figure}[H]
\centering
  \textbf{Dynamics of closed-loop system with output feedback }\par\medskip
\begin{center}
\subfloat[$w(x,t)$]{\includegraphics[scale=0.25]{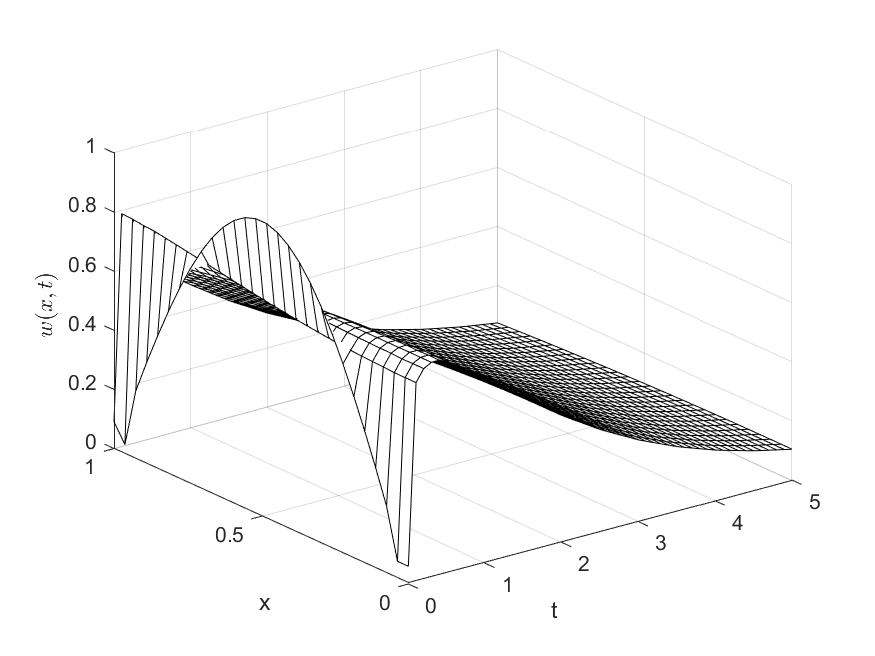}}
\hfill
\subfloat[$v(x,t)$]{\includegraphics[scale=0.25]{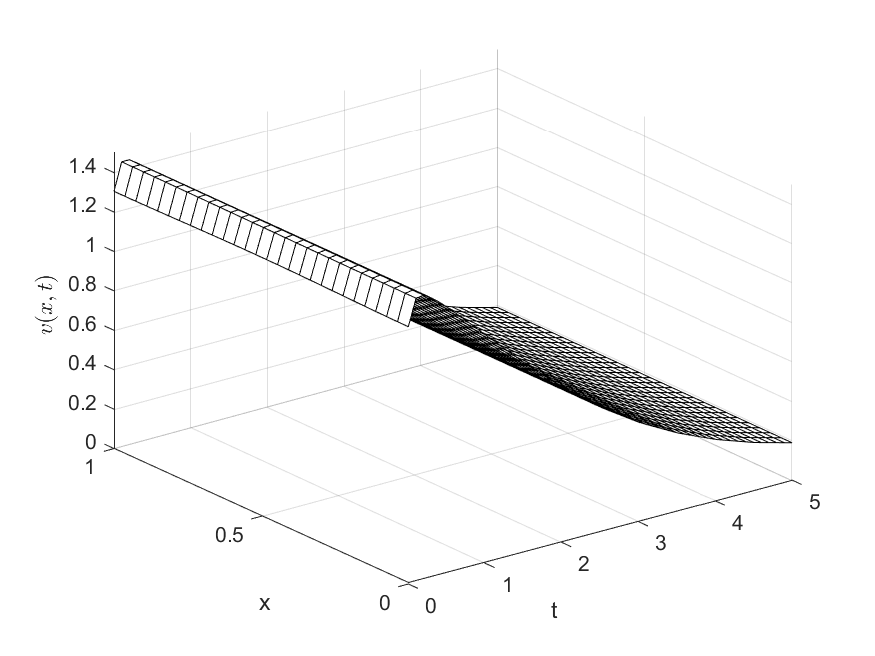}}
\end{center}
\caption{{A 3D landscape of the dynamics of a controlled coupled parabolic-elliptic system (\ref{1})-(\ref{4}) after applying the output feedback control input \eqref{cs}. The uncontrolled system is unstable, but the use of output feedback leads to an exponentially stable closed-loop system.}}
\label{fig-outputfeedback-1}
\end{figure}

\begin{figure}[H]
\begin{center}
 \textbf{Comparison of open-loop and closed-loop  $L_2$-norms with output feedback }\par\medskip
\subfloat[$L_2$-norm $\|w(\cdot,t)\|$]{\includegraphics[scale=0.25]{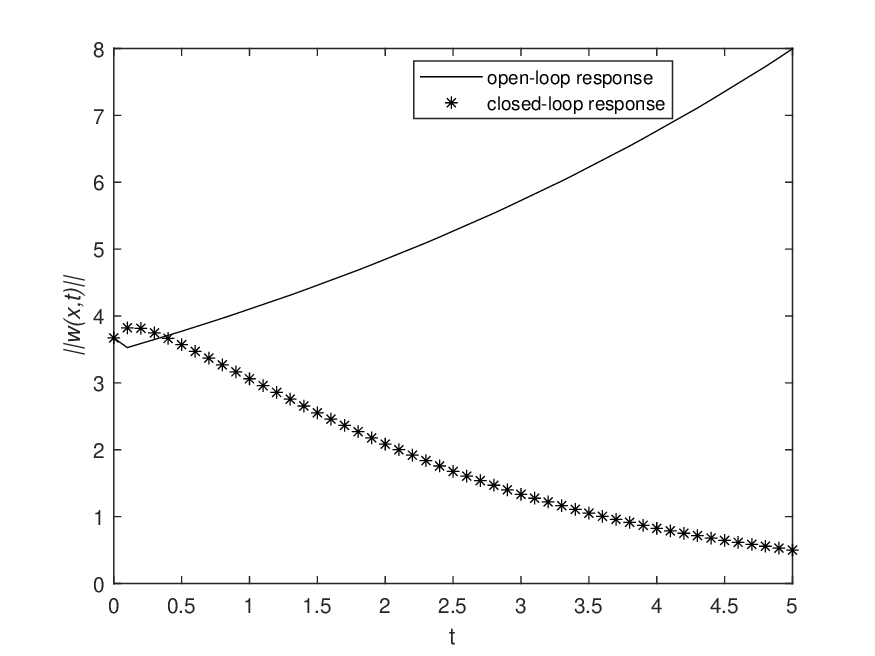}}
\hfill
\subfloat[$L_2$-norm  $\|v(\cdot,t)\|$]{\includegraphics[scale=0.25]{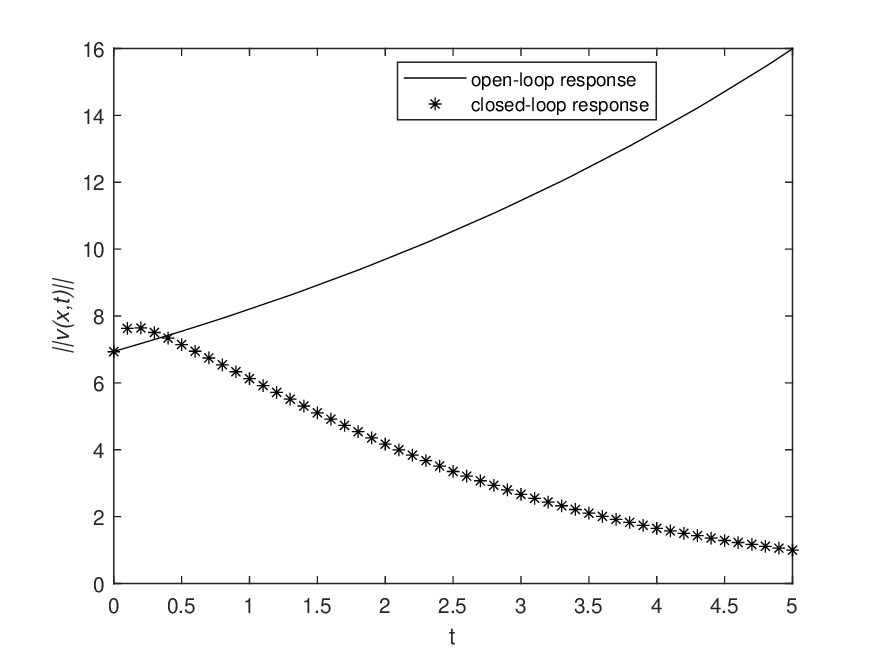} }
\end{center}
\caption{ A comparison between the $L_2$-norms of  $w(x,t)$ and $v(x,t)$ for the uncontrolled system and the system controlled with output feedback.  Without control, the norms of both states grow, due to instability. However, with output feedback control, both states decay to zero.}
\label{fig-outputfeedback-2}
\end{figure}

To achieve stability, the control gain  $c_2=1.2-\rho$, ensuring that the stability condition \eqref{bound on c1} is satisfied.  Additionally,  to ensure exponential stability of the  error dynamics $o_2=3$, which satisfies inequality \eqref{critera for c1}.  However, after applying control, the state of the coupled system converges to the steady-state solution; see Figure\ref{fig-outputfeedback-1}. This convergence is clearly depicted in the comparison shown in Figure\ref{fig-outputfeedback-2}, where we compare the $L_2$- norm of the controlled and uncontrolled states.

\section{Conclusion}\label{section5}
Stabilization of systems composed of coupled parabolic and elliptic equations presents considerable challenges. The first part of this paper considers  the boundary stabilization of a linear coupled parabolic-elliptic system.  Previous literature has shown that coupling between the equations can result in an unstable system.  Previous work on stabilization via two boundary control inputs was used  in \cite{krstic2008boundary}. 
In this paper we used  a single control input to stabilize both equations. Deriving one control law that stabilizes the system posed several challenges. One approach is   to rewrite  the coupled system into one equation in terms of the parabolic state. But due  to the appearance of a Fredholm operator, this makes it  difficult to establish a suitable kernel for the backstepping transformation.  Using  separate transformations for each  of the parabolic and the elliptic states would be  quite complicated. In this paper we transformed only the   parabolic part of the system, which simplified the calculations. This enabled reuse of  a previously calculated  transformation. However, the price is that  this transformation mapped the original coupled system into a complicated target system which further complicated showing stability of the target system.  Lyapunov theory was used to obtain a sufficient condition for stability of the target system. 

The second part of the paper focused on the observer design problem. Several synthesises are proposed depending on the available measurements. Output injections  were chosen so that the exponential stability of the observation error dynamics is ensured.  Again, instead of looking for a new state transformation that maps the original error dynamics into an exponentially stable target system,  well-known transformations in the literature are employed. Then, the exponential stability of the original error dynamics is shown by  establishing suitable sufficient conditions for stability.   The key to obtaining a  stability condition was again to  use Lyapunov theory. As for controller design, the technical conditions for  observer design  depend on the number of the available measurements. When  measurements for both states were provided, two transformations were applied to both parabolic and elliptic  states of the error dynamics. A total of four filters, two throughout the domain and two  at the boundary, were needed. On the other hand, when a single measurement for the parabolic state is given, one boundary filter was designed for the  parabolic equation. However, in the latter case, a more restrictive condition for stability of the error dynamics was obtained.
Observer design with a  single sensor parallels to a great extent that of  stabilization via one control signal. Showing exponential stability of the target error dynamics, in the situation when one measurement is provided,  results in a very restrictive constraint on the parameters of the coupled system. 
%However, determining a sufficient condition for stability by studying the original error dynamics, where the obtained output injections are incorporated, gives a more relaxed criterion. 

In a final section, control and observation were combined to obtain an output feedback controller. The results were again illustrated with simulations. 

 In this paper, the control as well as the observer gains were designed at $x=1$. The  situation is similar when the control and the observer gains are placed at $x=0$, although different backstepping transformations are needed. This is covered in detail in \cite{Ala}. 

Open questions include finding   weaker sufficient conditions for stability of both the controller and the observer.
Future work is aimed at improving the conditions on the control parameter $c_2$ concerning the boundary stabilization problem. Similarly, relaxing stability condition on the parameters for the observer design problem with a single measurement is  of interest.  In many problems, the equations are nonlinear and  design of a boundary control for a  nonlinear coupled parabolic-elliptic equations will also be studied.

\bibliographystyle{plain}        % Include this if you use bibtex 
\bibliography{autosam}           % and a bib file to produce the 
                                 % bibliography (preferred). The
                                 % correct style is generated by
                                 % Elsevier at the time of printing.

%\begin{thebibliography}{99}     % Otherwise use the  
                                 % thebibliography environment.
                                 % Insert the full references here.
                                 % See a recent issue of Automatica 
                                 % for the style.
%  \bibitem[Heritage, 1992]{Heritage:92}
%     (1992) {\it The American Heritage. 
%     Dictionary of the American Language.}
%     Houghton Mifflin Company.
%  \bibitem[Able, 1956]{Abl:56}
%     B.~C.~Able (1956). Nucleic acid content of macroscope. 
%     {\it Nature 2}, 7--9. 
%  \bibitem[Able {\em et al.}, 1954]{AbTaRu:54}   
%     B.~C. Able, R.~A. Tagg, and M.~Rush (1954).
%     Enzyme-catalyzed cellular transanimations.
%     In A.~F.~Round, editor, 
%     {\it Advances in Enzymology Vol. 2} (125--247). 
%     New York, Academic Press.
%  \bibitem[R.~Keohane, 1958]{Keo:58}
%     R.~Keohane (1958).
%     {\it Power and Interdependence: 
%     World Politics in Transition.}
%     Boston, Little, Brown \& Co.
%  \bibitem[Powers, 1985]{Pow:85}
%     T.~Powers (1985).
%     Is there a way out?
%     {\it Harpers, June 1985}, 35--47.

%\end{thebibliography}

\appendix
%\section{A summary of Latin grammar}    % Each appendix must have a short title.
%\section{Some Latin vocabulary}         % Sections and subsections are supported  
                                        % in the appendices.
\end{document}